%% file: _0FLM-ESOP.tex
\documentclass[11pt,a4paper,fleqn]{article}
\pdfoutput=1

\RequirePackage{amsthm,amsmath,natbib}
\RequirePackage{amssymb,amstext}
\RequirePackage{hypernat}
\usepackage[pdfpagemode=UseOutlines ,plainpages=false
,hypertexnames=false ,pdfpagelabels ,hyperindex=true,colorlinks=true]{hyperref}
	\makeatletter%
	\Hy@breaklinkstrue%
	\makeatother%
\usepackage{color}
\definecolor{darkred}{rgb}{0.6,0.0,0.1}
\definecolor{darkgreen}{rgb}{0,0.5,0}
\definecolor{darkblue}{rgb}{0,0,0.5}
\hypersetup{colorlinks ,linkcolor=darkblue ,filecolor=darkgreen
,urlcolor=darkblue ,citecolor=black}

\usepackage{url}
\usepackage{verbatim}
\usepackage{ bm}

\renewcommand{\cite}{\citet}
\usepackage{jan-abrev-package}

\definecolor{dgreen}{rgb}{0,0.5,0}
\definecolor{dblue}{rgb}{0,0,0.9}
\definecolor{dred}{rgb}{0.6,0.0,0.1}
\definecolor{dgold}{rgb}{0.5,0.3,0.0}
\definecolor{dvio}{rgb}{0.6,0.3,0.5}
\definecolor{gray}{rgb}{0.5,0.5,0.5}

\newtheoremstyle{mysc}% name
  {3pt}%      Space above
  {3pt}%      Space below
  {\it}%         Body font
  {}%         Indent amount (empty = no indent, \parindent = para indent)
  {\color{darkred}\sc}% Thm head font
  {.}%        Punctuation after thm head
  {.5em}%     Space after thm head: " " = normal interword space;
  {}%         Thm head spec (can be left empty, meaning `normal')

\newtheoremstyle{myex}% name
  {10pt}%      Space above
  {10pt}%      Space below
  {\rm}%         Body font
  {}%         Indent amount (empty = no indent, \parindent = para indent)
  {\color{darkred}\sc}% Thm head font
  {.}%        Punctuation after thm head
  {.5em}%     Space after thm head: " " = normal interword space;
  {}%         Thm head spec (can be left empty, meaning `normal')

\theoremstyle{mysc}\newtheorem{prop}{Proposition}[section]
\theoremstyle{mysc}\newtheorem{assumption}{Assumption}[section]
\theoremstyle{mysc}
\theoremstyle{mysc}\newtheorem{theo}[prop]{Theorem}
\theoremstyle{mysc}
\theoremstyle{mysc}\newtheorem{lem}[prop]{Lemma}
\theoremstyle{myex}\newtheorem{rem}{Remark}[section]
\theoremstyle{myex}
\theoremstyle{myex}

\numberwithin{equation}{section}

\author{{\sc Jan Johannes}\thanks{Institut f\"ur Angewandte Mathematik, Im Neuenheimer Feld, 294, D-69120
    Heidelberg, Germany, e-mail:
    \url{johannes@math.uni-heidelberg.de}}\\[.5ex]
  Ruprecht-Karls-Universit\"at Heidelberg}

\title{{\bf Functional linear instrumental regression under  second order stationarity.}}

\begin{document}
\date{~}
\maketitle

%\vskip 1cm
\begin{abstract} We consider the problem of estimating the slope parameter in  functional linear instrumental regression, where in the presence of an instrument $W$, i.e., an exogenous random function, a scalar response
$Y$ is modeled in dependence of an endogenous random function
$X$. Assuming second order stationarity jointly for $X$ and $W$  a nonparametric estimator of  the functional slope parameter and its derivatives 
is proposed based on an $n$-sample of $(Y,X,W)$. In this paper the minimax optimal rate of convergence of the estimator is derived  assuming that the slope parameter belongs to the well-known Sobolev space of periodic functions.  We discuss the cases that the cross-covariance operator associated to the random functions $X$ and $W$ is finitely, infinitely or in some general form smoothing.  
\end{abstract}

%\vspace{0.5cm}
\begin{tabbing}
\noindent \emph{Keywords:} \=Functional linear model, Instrument,  Orthogonal series estimation,  \\
\>Spectral cut-off, Optimal rate of convergence, Sobolev space.\\[.2ex]
\noindent\emph{JEL classifications:} Primary C14; secondary C30.
\end{tabbing}
The author gratefully acknowledges support from the research program "New Challenges for New Data" of LCL and Genes.
\input{_1intro}
\input{_2formal}
\input{_3theo-prop-gen}
\input{_4theo-prop-infin}
\input{_5concl}
\appendix
\section{Appendix: Proofs}\label{app:proofs}
\subsection{Proofs of Section \ref{sec:gen}}\label{app:proofs:gen}
\input{_a1proofs-consistence}
\input{_a2proofs-lower}
\input{_a3proofs-upper}
\input{_a4proofs-technical}
\subsection{Proofs of Section \ref{sec:infin}}\label{app:proofs:infin}
\input{_a5proofs-infin}

\bibliography{FLM-ESOP}
\end{document}

%% file: _1intro.tex
\section{Introduction}\label{sec:intro}
The analysis of functional data is becoming very important in a diverse range of disciplines,
including medicine, linguistics, chemometrics as well as econometrics (see for instance 
\cite{RamsaySilverman2005} and \cite{FerratyVieu2006}, for several case studies).  In particular, there is a wide diversity of applications in economics. 
\cite{ForniReichlin1998} study business cycle dynamics and \cite{PredaSaporta2005}
 consider shares at the Paris stock exchange, to name but a few. Roughly speaking, in all
these applications the dependence 
of a response variable $Y$  on the variation of an explanatory random function $X$ is modeled by a functional linear regression model, that is,
\begin{equation}\label{intro:e1}Y=\int_{0}^1\beta(t)X(t)dt +\sigma U,\quad\sigma>0,
\end{equation}
for some error term $U$. The  important point to note here is  that often in economical applications the commonly used hypothesis,  that the  regressor $X$ is exogenous, can  be rejected using, for example,  a  test  proposed by \cite{BH07res}. Thus analyzing the influence of endogeneity  is  of particular interest in econometrics.   One objective   is then to estimate nonparametrically in the presence of an instrument $W$ the slope function $\beta$ or its derivatives  based on an $n$-sample of $(Y,X,W)$.  
\paragraph{Background.} Suppose first the regressor $X$ is exogenous, i.e., $\Ex[UX(s)]=0$, $s\in[0,1]$. In this case the estimation of the slope function $\beta$ has been considered, for example,  in \cite{CardotFerratySarda2003}, \cite{MullerStadtmuller2005},       
 \cite{HallHorowitz2007} or \cite{CrambesKneipSarda2007}. Assuming  the  random function $X$ to be centered the most popular approach is to  multiply both sides in (\ref{intro:e1}) by $X(s)$.  Then taking the expectation leads to 
\begin{equation}\label{intro:e2}\Ex[YX(s)]=\int_{0}^1\beta(t)\Cov(X(t),X(s))dt,\quad s\in[0,1].
\end{equation}
 The normal equation (\ref{intro:e2}) is the continuous equivalent of a normal equation in a classical  linear model. To be more precise, suppose the random function $X$ and the slope function $\beta$ to be square integrable, then  their generalized Fourier coefficients $X_i:=\int_0^1 X(s)\psi(s)ds$ and $\beta_i:=\int_0^1 \beta(s)\psi(s)ds$, $i\in\N$, with respect to some  orthonormal basis $\{\psi_i\}$ are well-defined. The functional linear model \eqref{intro:e1} (FLM for short) and hence the normal equation (\ref{intro:e2}) can be rewritten as 
\begin{equation}\label{intro:e3}Y=\sum_{i=1}^\infty \beta_i X_i+\sigma U,\quad\mbox{ and } \quad
\Ex[YX_j] = \sum_{i=1}^\infty \beta_i \cdot \Cov(X_i,X_j),\quad j=1,2,\dotsc,
\end{equation}
respectively. Therefore, the  FLM \eqref{intro:e1} extends the  linear model  (LM for short) $Y=\sum_{i=1}^k \beta_i X_i+\sigma U$, $k\in\N$, to an infinite number of regressors. Since in analogy to  the estimation in  the LM 
recovering from \eqref{intro:e3}  the coefficients $(\beta_i)_{i\in\N}$  
necessitates the inversion of the infinite dimensional covariance matrix $\Sigma_\infty:=(\Cov(X_i,X_j))_{i,j\in\N}$,  the estimation of $\beta$ is called an inverse problem. It is well-known that in both, the linear and the functional linear  model  identification as well as the accuracy of any estimator  depends strongly on the properties of the covariance matrix $\Sigma_k:=(\Cov(X_i,X_j))_{i,j=1}^k$ and $\Sigma_\infty$ respectively. That is, in both cases the coefficients  can be identified as long as  the covariance matrix $\Sigma_k$ and $\Sigma_\infty$ respectively, is not singular. Moreover, in the LM  a high degree of multicolinearity  between the regressors $X_1,\dotsc,X_k$, that is, the smallest  eigenvalue of $\Sigma_k$ is close to zero, produces unacceptable uncertainty in the coefficient estimates. However, as long as the covariance matrix $\Sigma_k$ is not singular an ordinary least squares estimator (LSE for short) will be consistent and leads under fairly weak assumptions to a minimal asymptotic variance.  In general the situation in the FLM is  different. Since under very mild assumptions zero is an accumulation point of the eigenvalues of $\Sigma_\infty$ we always have to face  a classical multicolinearity problem in the presence of many regressors. Therefore, although the covariance matrix $\Sigma_\infty$ is not singular the LSE will not longer  be   consistent.  This corresponds to the setup of ill-posed inverse problems. 

 In the literature several approaches are proposed in order to circumvent in the FLM the instability issue.  
Essentially, all of them replace  the covariance matrix $\Sigma_\infty$ in equation (\ref{intro:e3}) by a regularized version. 
A popular  example is based on the functional principal
components  regression (c.f. \cite{Bosq2000}, \cite{MullerStadtmuller2005} or  \cite{CardotMasSarda2007}), which corresponds to a  method  called spectral cut-off in the literature of numerical analysis  (c.f. \cite{Tautenhahn96}). Another example is the Tikhonov regularization (c.f. \cite{HallHorowitz2007}), where the regularized solution $\beta_\alpha$ is defined as unique minimizer of the Tikhonov functional 
$F_\alpha(\beta)=\sum_{j=1}^\infty\{\Ex[YX_j]-\sum_{i=1}^\infty \beta_i \cdot \Cov(X_i,X_j) \}^2 + \alpha \sum_{j=1}^\infty\beta_j^2$ for some strictly positive $\alpha$.  Regularization using
a penalized least squares approach after projection onto some basis (such as
splines) is also considered in \cite{RamsayDalzell1991}, \cite{EilersMarx1996} or \cite{CardotFerratySarda2003}.  The common aspect of all these regularization schemes is the introduction of an additional  regularization parameter $\alpha$ (for example, the parameter determining the weight of the penalty in the  Tikhonov functional). The risk of the resulting regularized  estimator   can then be decomposed, roughly speaking,  into a function of the risk of the  estimators of $\Ex[YX_j]$ and $\Cov(X_i,X_j)$, $i,j\in\N$, plus an additional bias term which is a function of the regularization parameter $\alpha$. The optimal value of $\alpha$ is then obtained by balancing these two terms. However, in order to obtain a rate of convergence 
additional regularity assumptions on the slope function $\beta$ and the  infinite dimensional covariance matrix $\Sigma_\infty:=(\Cov(X_i,X_j))_{i,j\in\N}$ are necessary (a detailed discussion in the context of inverse problems in econometrics can be found in \cite{CFR05handbook} or \cite{JohannesVan-BellegemVanhems2011}).

The objective of this paper is to study the estimation of the slope function $\beta$ when the
regressor $X$ is endogenous, which to the best of our knowledge has  not yet been considered in the literature. In the following the approach of this paper is described in more details. 
\paragraph{Methodology.} To treat the endogeneity problem, we assume that an instrument $W$, i.e., an exogenous random function, is given.   Assuming  the  random function $W$ to be centered and square integrable, 
we consider its generalized Fourier coefficients $W_i:=\int W(s)\psi'_i(s)ds$, $i\in\N$, with respect to some  orthonormal basis $\{\psi'_i\}$ not necessarily the same as $\{\psi_i\}$ used above in the decomposition of $X$ and $\beta$. Then multiplying the equation 
\eqref{intro:e1} by $W_j$  and taking the expectation leads to the normal equation
\begin{equation}\label{intro:e4}\Ex[YW_j] = \sum_{i=1}^\infty \beta_i \cdot \Cov(X_i,W_j),\quad j=1,2,\dotsc.,
\end{equation}
which  provides a natural extension of the linear instrumental regression (LIR for short) $Y=\sum_{i=1}^k \beta_i X_i+\sigma U$  with $\Ex[W_jU]=0$, $j=1,\dotsc,q$,  to an infinite number of regressors and instruments. Therefore, in the presence of an instrument $W$ we call \eqref{intro:e1}  functional linear instrumental regression (FLIR for short). The estimation of the  coefficients in both,  linear and  functional linear  instrumental regression is  then again an inverse problem, since 
it involves now the inversion of the cross-covariance matrix $\Sigma_{kq}:=(\Cov(X_i,W_j))_{i,j=1}^{k,q}$ and $\Sigma_{\infty\infty}:=(\Cov(X_i,W_j))_{i,j\in\N}$ respectively.
 Furthermore, in both cases the coefficients  are identifiable as long
 as   $\Sigma_{kq}$ and $\Sigma_{\infty\infty}$ respectively, is not
 singular and moreover, the obtainable accuracy of any estimator
 depends now on the properties  $\Sigma_{kq}$ and
 $\Sigma_{\infty\infty}$ respectively. It is worth to pointing out
 that the FLIR parallels  developments in econometric theory
 such as   nonparametric  instrumental regression
 (c.f. \cite{DarollesFanFlorensRenault2011}, \cite{NP03econometrica},
 \cite{HH05aos} or \cite{FlorensJohannesVan-Bellegem2010}),   nonparametric instrumental quantile
regression of \cite{HL07econometrica} or semi-nonparametric estimation of Engel curve
with shape-invariant specification of \cite{BCK07econometrica}.

The estimator of  the slope function in FLIR considered in this paper is based on a two stage least squares approach. To be more precise, consider first the LIR. Then as long as the cross-covariance matrix $\Sigma_{kq}$ is not singular a two stage least squares procedure (2SLS for short) will lead to a consistent estimator. That is, in a first step a linear regression of the endogenous vector  $X=(X_1,\dotsc,X_k)^t$ onto the vector of instruments $W=(W_1,\dotsc,W_q)^t$ is performed, resulting into an estimator $\widehat{W}$ of the optimal linear instrument $\widetilde{W}$, i.e.,  the best linear predictor $\widetilde{W}:=\Sigma_{kq}\Sigma_q^{-1}W$ of $X$ with $\Sigma_{q}:=(\Cov(W_i,W_j))_{i,j=1}^{q}$. Note that  the optimal linear instrument  is well-defined as long as the covariance matrix $\Sigma_q$ of $W$ has full rank.
 Then in the second step an estimator of the  $k$-vector of coefficients $(\beta_j)$ is obtain considering a linear regression of $Y$ onto  $\widehat{W}$.  Applying a 2SLS approach in FLIR we have to face  additional technical difficulties given through the facts that  the optimal linear instrument, i.e., the best linear predictor of the random function $X$ given the random function $W$, is not always well-defined and that both stages of the estimation procedure necessitate the solution of  an ill-posed inverse problem (see the discussion above in case of an exogenous regressor). Therefore, assuming the optimal linear instrument is well-defined, we apply in each stage a regularization scheme in order to circumvent the instability issue.   Although the estimation in the first step has to be  stabilized, it has only a minor influence on the obtainable accuracy of the final estimator. Particularly,  the proposed estimator of $\widetilde{W}$ will in general not be  optimal.  The main complexity of the estimation problem is contained in the second stage. 
To be more precise, if the optimal linear instrument $\widetilde{W}$ is given, then the second step in fact only consists of the  estimation  in a FLM  \eqref{intro:e1} given now with exogenous regressor $\widetilde{W}$. Thereby,  the relationship between the regularity assumption on the slope function $\beta$ and the infinite dimensional covariance matrix $\widetilde{\Sigma}_\infty:=(\Cov(\widetilde{W}_i,\widetilde{W}_j))_{i,j\in\N}$   associated to the  instrument $\widetilde{W}$ determines the obtainable accuracy of any estimator of $\beta$ (see also the discussion above in case of an exogenous regressor). Nevertheless, the  instrument $\widetilde{W}$ is not given and thus has to be estimated. However,
 the estimation in the first step is possible without changing the optimal rate of the estimator of $\beta$, where only higher moment conditions are the price to pay.

Suppose the slope function $\beta$ belongs to the Sobolev space  of periodic functions $\cW_p$ (defined below). Given  an $n$-sample of $(Y,X,W)$ our objective is not only the estimation of the slope function $\beta$ itself but also of its derivatives. We show that the relationship between the Sobolev spaces and the covariance matrix $\widetilde{\Sigma}_\infty$ associated to the optimal linear instrument $\widetilde{W}$, i.e., the \lq\lq smoothing\rq\rq\ property of $\widetilde{\Sigma}_\infty$, is essentially determining  the optimal rate of convergence of any estimator.   We now describe two examples. First consider the covariance matrix $\widetilde{\Sigma}_\infty$  to be finitely smoothing, that is, the range of  $\widetilde{\Sigma}_\infty$ equals $\cW_a$ for some $a>0$. Then the optimal rate is a polynomial of the sample size $n$. 
It is worth to note that all published results in the FLM with exogenous regressor consider only this case (c.f. \cite{HallHorowitz2007} and  \cite{CrambesKneipSarda2007}). However,  assuming  $\widetilde{\Sigma}_\infty$ to be finitely smoothing excludes several interesting situations, such as our second  example. Suppose $\widetilde{\Sigma}_\infty$ to be   infinitely smoothing, that is, the range of $|\log(\widetilde{\Sigma}_\infty)|^{-1}$ equals $\cW_a$ for some $a>0$. Then the optimal rate is a logarithm of the sample size $n$. The important point to note here is the theory behind these cases can be generalized by  using an index function $\kappa$ (c.f. \cite{NairPereverzevTautenhahn05}), which `links' the range of $\widetilde{\Sigma}_\infty$ and the Sobolev spaces. Then $\widetilde{\Sigma}_\infty$ is called in some general form smoothing and moreover the index function $\kappa$ determines the functional form of the optimal rate of convergence.  A similar approach in the context of nonparametric instrumental regression, where the conditional expectation  plays the same role as the covariance matrix $\widetilde{\Sigma}_\infty$, can be found  in  \cite{ChenReiss2011} or \cite{JohannesVan-BellegemVanhems2011}.

In this paper we deal with the estimation of the slope function when the regressor $X$ and the instrument $W$ are jointly second order stationary (defined below). We derive a lower bound of the rate of convergence for any estimator of $\beta$ or its derivatives assuming some general form of smoothing of $\widetilde{\Sigma}_\infty$.  Assuming second order stationarity we propose an orthogonal series estimator of $\beta$ and its derivatives based on a spectral cut-off (thresholding in the Frequency domain). Then we show that the rate of the lower bound provides also an upper bound for the risk of  the orthogonal series estimator. Therefore, the rate is optimal and hence the proposed estimator is minimax-optimal. The results for  general smoothing $\widetilde{\Sigma}_\infty$  imply then as propositions the  minimax optimal rate of convergence in estimating $\beta$ and its derivatives respectively in case of  finitely  as well as infinitely smoothing $\widetilde{\Sigma}_\infty$.

\paragraph{Organization of the paper.} We summarize  in Section 2 the
model assumptions and define the estimator of $\beta$ and its
derivatives. In Section 3 we provide minimal conditions to ensure
consistency of the estimator. Furthermore, we derive a lower and an
upper bound for the risk in the Sobolev norm when
$\widetilde{\Sigma}_\infty$ is in some general form smoothing. This
results are illustrated in Section 4 assuming
$\widetilde{\Sigma}_\infty$ to be finitely or infinitely smoothing and
Section 5 concludes. All proofs can be found in the Appendix. 

%% file: _2formal.tex
\section{Formalization of the model  and definition of the estimator}\label{sec:methodology}
\paragraph{Model.}The setting of this paper can be summarized through the model
      \begin{subequations}
	\begin{equation}
	\label{model}
	Y = \int_0^1\beta(t) X(t)dt + \sigma U, \quad \sigma>0,
	\end{equation}
where  $Y \in \R$ is a response variable, the endogenous random function $X$ is defined on the interval $[0,1]$ and $U$ is a centered error term with variance one such that 
	\begin{align}
	\label{model:W}
	\Ex [ U W(t) ] = 0,\quad t\in[0,1]
	\end{align}
      \end{subequations}
for some instrument $W$, i.e., an exogenous random function defined also on $[0,1]$. The objective is the nonparametric estimation of the
slope function $\beta$ and its derivatives based on a $n$-sample of $(Y,X,W)$. We assume throughout the paper that the random functions $X$ and $W$ are defined on the interval $[0,1]$ that (technically) simplifies the notations. Of course, it does not touch the applicability of the model and suggested estimator in a general setting when $X$ and $W$ are defined on some compact intervals $I_1$ and $I_2$, respectively. Moreover, we suppose  that the random functions $X$ and $W$  have a finite second moment, i.e., $\int_0^1\Ex|X(t)|^2dt<\infty$ and $\int_0^1\Ex |W(t)|^2dt<\infty$.  In order to simplify the presentation we assume that the mean function of $X$ and $W$ are zero. Then multiplying both sides in (\ref{model}) by $W(s)$, $s\in[0,1]$, and taking the expectation leads to 
\begin{equation}\label{cross-cov}\Ex[YW(s)]=\int_{0}^1\beta(t)\Cov[X(t),W(s)]dt=:[T_{_{WX}} \beta](s),\quad s\in[0,1],
\end{equation}
where the function $\Ex[YW(\cdot)]$ is square integrable and $T_{_{WX}}$ denotes the cross-covariance operator associated to the random functions $X$ and $W$. Note that the cross-covariance matrix $\Sigma_{\infty\infty}$ considered in the introduction satisfies $\Sigma_{\infty\infty}=(\int_0^1 [T_{_{WX}}\psi_i](s)\psi_j'(s)ds)_{i,j\in\N}$. Estimation of $\beta$ is thus linked with the inversion of the cross-covariance operator $T_{_{WX}}$ of $(X,W)$, and, hence called an inverse problem.  Throughout the paper we require the following assumption, which  
 provides   a   necessary and sufficient condition for the existence of a unique solution of equation (\ref{cross-cov}).
\begin{assumption}\label{ass:ident} The cross-covariance operator $T_{_{WX}}$ associated to the random functions $X$ and $W$ is injective and the function $\Ex[YW(\cdot)]$ belongs to the range $\cR(T_{_{WX}})$ of $T_{_{WX}}$.\end{assumption} 
In case  a solution of the normal equation (\ref{cross-cov}) does not exist all the results below can also straightforward be obtained for the unique  least-square solution with minimal norm, which exists if only  if $\Ex[YW(\cdot)]$ is contained in the direct sum of $\cR(T_{_{WX}})$ and its orthogonal complement $\cR(T_{_{WX}})^\perp$ (for a definition and detailed discussion in the context of inverse problems c.f. \cite{EHN00} or \cite{CFR05handbook}).

\paragraph{Notations and basic assumptions.}In this paper we suppose that  the  random function $(X,W)$ is second order stationary and, hence there exists a function $c_{_{WX}}:[-1,1]\to\R$ such that  $\Cov[X(t),W(s)]=c_{_{WX}}(s-t)$, $t,s\in[0,1]$. Notice that due to the finite second moment of $X$ and $W$ the cross-covariance function $c_{_{WX}}(\cdot)$ is square integrable. Therefore, its Fourier coefficients with respect to the  Fourier complex exponentials, i.e., 
\begin{equation}\label{cov:sv}
c_k:= \int_{-1}^1 c_{_{WX}}(t) \exp(-2\pi k it) dt,\quad\mbox{ for all }k\in\Z,\end{equation}
are well-defined and by applying the well-known convolution theorem we have
\begin{equation}\label{cov:svd}T_{_{WX}}\phi_k=c_k \cdot\phi_k\quad\mbox{with}\quad \phi_k(t):= \exp(2\pi k it),\; t\in[0,1],\quad\mbox{ for all }k\in\Z.\end{equation}  
Thereby, it is convenient to consider the real-valued  random functions  $X$ and $W$ as elements of the Hilbert space $L^2[0,1]$ of square integrable complex valued functions defined on $[0,1]$, which is endowed with  inner product $\skalarV{f,g}=\int_0^1f(t)\overline{g(t)}dt$  and associated norm $\normV{f}=\skalarV{f,f}^{1/2},$ $f,g\in L^2[0,1]$. Here and subsequently, $\overline{g(t)}$ denotes the complex conjugate of $g(t)$. Furthermore, the cross-covariance operator $T_{_{WX}}$ is a well-defined mapping from  $L^2[0,1]$ into itself.
Consider the centered  complex valued random variables $\skalarV{\phi_k,X}$ and $\skalarV{W,\phi_k}$, $k\in\Z$, which due to the identity \eqref{cov:svd} satisfy%
\begin{equation}
 c_k= \Ex[\skalarV{\phi_k,X}\skalarV{W,\phi_k}]\quad\mbox{ and }\quad 0= \Ex[\skalarV{\phi_k,X}\skalarV{W,\phi_j}]\quad \mbox{ for all }j\ne k.
\end{equation}
Now an equivalent formulation of Assumption \ref{ass:ident} is given by 
\begin{equation}\label{ass:ident:e1}
|c_k|^2>0,\quad\mbox{for all }k\in\Z\qquad  \mbox{ and }\quad \sum_{k\in\Z} \frac{|\Ex[Y\skalarV{W,\phi_k}]|^2}{ |c_k|^2}<\infty. 
\end{equation}
 \paragraph{Optimal linear instrument.}Let $x_k:= \Var\skalarV{X,\phi_k}$, $w_k:= \Var\skalarV{W,\phi_k}$  and define $\lambda_k:= c_k^2 / w_k\leq x_k$, $k\in \Z$, where due to the finite second moment of $X$, i.e., $\Ex\normV{X}^2=\sum_{k\in\Z}x_k$, the sequences $(\lambda_k)_{k\in\Z}$ is summable. If we further assume  that  $\sup_{k\in\Z}|\lambda_k/ w_k|<\infty$, then the complex valued random function 
\begin{equation}\label{bm:oiv}\widetilde{W}:=\ell(W):=\sum_{k\in\Z} \frac{\overline{c}_k}{w_k}\cdot \skalarV{W,\phi_k}\cdot\phi_k\end{equation} is well-defined, i.e., $\normV{\widetilde{W}}<\infty$. Note that $\ell$ is a linear operator mapping $L^2[0,1]$ into itself. If in addition $\sum_{k\in\Z}\lambda_k/ w_k<\infty$, then  $\ell$ is a Hilbert-Schmidt operator and $\widetilde{W}=\ell(W)$ is the best linear predictor of $X$ based on $W$. That is, $\ell$ minimizes the mean prediction error $\Ex\normV{X-\ell'(W)}^2$ over all   Hilbert-Schmidt operator $\ell'$ (c.f. \cite{Bosq2000}). Therefore, we call $\widetilde{W}$  optimal linear instrument.   Throughout the paper we  suppose the  linear predictor $\widetilde{W}$ is well-defined, i.e., $\sup_{k\in\Z}|\lambda_k/ w_k|<\infty$, which implies an additional restriction on the behavior of the sequences $(\lambda_k)_{k\in\Z}$  and $(w_k)_{k\in\Z}$ as $|k|\to\infty$. In particular  the sequence of variances $(w_k)_{k\in\Z}$  associated to the instrument $W$ has to tend slowlier to zero than the  sequence $(c_k)_{k\in\Z}$ of cross-covariances   associated to $X$ and  $W$. Note that although we suppose the optimal linear instrument $\widetilde{W}$ exist, in general it is not known to the econometrician.
%%%%%%%%%%%%%%%%
\paragraph{Moment assumptions.}The results  derived below involve  additional  conditions on the moments of the  random functions $X$ and $W$ and the error term $U$, which we formalize now. Let $\cF$ be the set of all centered  second order stationary random functions $(X,W)$.
Here and subsequently,  $\cF^{m}_{\eta,\tau}$, $m\in\N$, $\eta,\tau\geq1$, denotes the subset of $\cF$  containing all random functions $(X,W)$ such that the $m$-th moment of the corresponding random variables  $\{\skalarV{X,\phi_k}/\sqrt{x_k}\}$ and $\{\skalarV{W,\phi_k}/\sqrt{w_k}\}$  are  uniformly bounded  and such that the  linear predictor of $X$ based on $W$ is well-defined, that is
\begin{multline}\label{bm:def:F}
\cF^{m}_{\eta,\tau}:=\Bigl\{ (X,W)\in\cF\;\text{ with }\sup_{k\in\Z} \Ex\Bigl|\frac{\skalarV{X,\phi_k}}{\sqrt{x_k}}\Bigr|^m \leq \eta\mbox{ and }\sup_{k\in\Z} \Ex\Bigl|\frac{\skalarV{W,\phi_k}}{\sqrt{w_k}}\Bigr|^m \leq \eta\\
\text{ and  associated values }(\lambda_k)_{k\in\Z} \text{ such that } 1\vee\sup_{k\in\Z}|\lambda_k/ w_k|\leq \tau\Bigr\}.
\end{multline}
In what follows, $\cE^m_\eta$ stands for the set of all centered error terms $U$ with variance one and finite $m$-th moment, i.e., $\Ex|U|^m\leq \eta$.
%%%%%%%%%%%%%%%%%
\paragraph{Estimation of $\beta$ as an ill-posed inverse problem.}  Consider the optimal linear instrument   $\widetilde{W}$ defined in \eqref{bm:oiv}, 
then due to  Assumption \ref{ass:ident} the normal equation \eqref{cross-cov} implies%
\begin{equation}\label{bm:e1}
\beta = \sum_{k\in\Z} \frac{g_k}{\lambda_k}\cdot \phi_k\quad\mbox{ with }g_k:=\skalarV{g,\phi_k},\; k\in\Z,\;\mbox{ and } g:=\Ex [Y\widetilde{W}(\cdot)].\end{equation}
Moreover,  $\lambda_k=\Ex|\skalarV{\widetilde{W},\phi_k}|^2,$ $k\in\Z$, are the eigenvalues of the covariance operator $T_{\widetilde{W}}$ associated to  $\widetilde{W}$. In other words,  the estimation of $\beta$ necessitates the inversion of the covariance operator $T_{\widetilde{W}}$. Accordingly,   replacing in \eqref{bm:e1} the unknown function $g$ by a consistent estimator $\widehat{g}$ does in general not lead  to a consistent estimator  of $\beta$  even in case of known  values $\{\lambda_k\}$. To be more precise, since the sequence $(\lambda_k)_{k\in\Z}$ tends to zero as $|k|\to\infty$,
$\Ex\normV{\widehat{g}- g}^2=o(1)$ does generally not imply $\sum_{k\in\Z}^\infty |\lambda_k|^{-2}\cdot\Ex|\skalarV{\widehat{g}-g,\phi_k}|^2 =o(1)$. Consequently, the estimation in  FLIR is called ill-posed and  additional regularity assumptions on the slope function $\beta$ are necessary in order to obtain a uniform rate of convergence (c.f. \cite{EHN00}). 

In this paper we assume that the slope function $\beta$ belongs  to the well-known  Sobolev space $\cW_{p}$, $p>0$, of periodic functions, which can be defined for $\nu\in\R$ by  
\begin{equation}\label{bm:def:Sob}
 \cW_{\nu}:=\Bigl\{f\in L^2[0,1]:\normV{f}_{\nu}^2 :=  \sum_{k\in\Z} \gamma_{k}^{\nu} |\skalarV{f,\phi_{k}}|^2 <\infty\Bigr\},
\end{equation} 
where $\{\phi_k\}$ are the complex exponentials  given in \eqref{cov:svd} and the  weights $\{\gamma_k\}$ satisfy%
\begin{equation}\label{bm:def:gamma}
\quad \gamma_{k}=1+|2\pi k|^2,\qquad k\in\Z.\end{equation}
Let $\cW_\nu^\rho:=\{f\in \cW_\nu:\normV{f}_\nu^2\leq\rho\}$ for $\rho>0$. Notice that for integer $\nu\in\N$ the  Sobolev space of periodic functions $\cW_\nu$ is equivalently given by  \begin{equation*}
 \cW_{\nu}=\Bigl\{f\in H_{\nu}: f^{(j)}(0)=f^{(j)}(1),\quad j=0,1,\dotsc,\nu-1\Bigr\},
 \end{equation*}
 where  $H_{\nu}:= \{ f\in L^2[0,1]:  f^{(\nu-1)}\mbox{ absolutely continuous }, f^{(\nu)}\in L^2[0,1]\}$  is a Sobolev space (c.f. \cite{Neubauer1988,Neubauer88}, \cite{MairRuymgaart96} or \cite{Tsybakov04}).

In the literature several approaches are proposed in order to circumvent an instability
issue due to the inversion of the covariance operator (for a detailed discussion in the context of inverse problems in econometrics we refer e.g. to \cite{CFR05handbook} and \cite{JohannesVan-BellegemVanhems2011}). Essentially, all of them replace  equation (\ref{bm:e1}) by a regularized version which avoids that the denominator becomes too small. For example, \cite{HallHorowitz2007} use in a functional linear model with exogenous regressor a Tikhonov regularization. There is a large number of alternative regularization schemes in the numerical analysis literature available like
the iterative Tikhonov regularization, Landweber iteration or the $\nu$-Method to name but a few (c.f. \cite{EHN00}). However, 
in this paper we regularize  equation (\ref{bm:e1}) by introducing a threshold $\alpha>0$ and  weights $\{\gamma_k\}$ defined in (\ref{bm:def:gamma}). For $\nu\geq 0$ we consider the regularized version ${\beta}_{\nu}$ given by
\begin{equation}\label{bm:def:sol:reg}
{\beta}_{\nu}:=\sum_{k\in\Z} \frac{{g}_{k}}{{\lambda}_{k}} \cdot\1\{{\lambda}_{k}/  \gamma_{k}^{\nu}\geq \alpha\}\cdot \phi_{k},\end{equation} 
which obviously  belongs  to the Sobolev space $\cW_{\nu}$. Thresholding in the Fourier domain in this situation is new, however has been used, for example, in a deconvolution problem in  \cite{MairRuymgaart96}, \cite{Neumann1997} or \cite{Johannes2009} and
coincides with an approach called spectral cut-off in the numerical analysis literature  (c.f. \cite{Tautenhahn96}). 
%%%%%%%%%%%%%%%
\paragraph{Definition of the estimator.}Let $(Y_1,X_1,W_1),\dotsc,(Y_n,X_n,W_n)$ be an   i.i.d. sample   of $(Y,X,W)$, which  we use in a first step to construct  an estimator $\widehat{W}_i$ of  the optimal linear instrument $\widetilde{W}_i$, $i=1,\dotsc,n$, exploiting the identity  \eqref{bm:oiv}. Consider the unbiased estimator of $c_{k}=\Ex \skalarV{\phi_{k},X}\skalarV{W,\phi_{k}}$ and   $w_{k}=\Ex |\skalarV{W,\phi_{k}}|^2$ given by
\begin{equation}\label{bm:def:est:OIV:c-w}
\widehat{c}_{k}:=\frac{1}{n}\sum_{i=1}^n \skalarV{\phi_{k},{X}_{i}}\skalarV{{W}_{i},\phi_{k}}\quad\mbox{ and }\quad\widehat{w}_{k}:=\frac{1}{n}\sum_{i=1}^n |\skalarV{{W}_{i},\phi_{k}}|^2,\quad k\in\Z,\end{equation}
respectively.  Then we define  the estimator of $\widetilde{W}_i$ by regularizing equation \eqref{bm:oiv}, i.e., by introducing a threshold $\alpha>0$, that is
\begin{equation}\label{bm:def:est:OIV}
\widehat{W}_{i}:=\sum_{k\in\Z}\frac{\overline{\widehat{c}_{k}}}{\widehat{w}_{k}}\cdot\1\{\widehat{w}_{k}\geq \alpha\}\cdot\skalarV{W_i,\phi_k}\cdot \phi_{k},\quad i=1,\cdots,n,
\end{equation}
 where the threshold $\alpha=\alpha(n)$ has to tend to zero as the sample size $n$  increases. In a second step we use the estimated optimal linear instrument to construct an estimator of $\beta$ based on the decomposition \eqref{bm:e1}. Consider the identities $\lambda_{k}=\Ex |\skalarV{\widetilde{W},\phi_{k}}|^2$ and $g_k=\Ex [Y\skalarV{\widetilde{W},\phi_k}]$, which motivate the estimators  defined by
\begin{equation}\label{bm:def:est:coeff:unknown}
\widehat{\lambda}_{k}:=\frac{1}{n}\sum_{i=1}^n |\skalarV{\widehat{W}_{i},\phi_{k}}|^2\quad\mbox{and}\quad
\widehat{g}_k:=\frac{1}{n}\sum_{i=1}^n Y_{i}\cdot\skalarV{\widehat{W}_i,\phi_k},\quad k\in\Z.\end{equation}
Finally, the estimator $\widehat{\beta}_\nu$ of $\beta$ is based on the regularized version \eqref{bm:def:sol:reg}. That is,% 
\begin{equation}\label{bm:def:est:reg}
\widehat{\beta}_{\nu}:=\sum_{k\in\Z}\frac{\widehat{g}_{k}}{\widehat{\lambda}_{k}} \cdot\1\{\widehat{\lambda}_{k}/  \gamma_{k}^{\nu}\geq \alpha\}\cdot \phi_{k},\end{equation}
which obviously  belongs also to the Sobolev space $\cW_{\nu}$.  
It is worth pointing out that  due to Parseval's formula  $\sum_{k\in\Z} \widehat{w}_{k} = \frac{1}{n} \sum_{i=1}^n \normV{{W}_i}^2$ is  finite. Thereby, the sum in (\ref{bm:def:est:OIV}) 
contains only a finite but random number of nonzero summands, and hence only a finite number of values $\widehat{\lambda}_{k}$ are nonzero. Consequently,  the sum in \eqref{bm:def:est:reg} contains only a finite number of nonzero summands.  We emphasize that the same threshold is used in the definition of  $\widehat{W}_i$ and $\widehat{\beta}_\nu$ given in \eqref{bm:def:est:OIV} and \eqref{bm:def:est:reg} respectively. In general this will not lead to an optimal  estimator  of $\widetilde{W}_i$, however as we will see below, it is sufficient to ensure the optimality of $\widehat{\beta}_\nu$.%

%% file: _3theo-prop-gen.tex
\section{Optimal estimation of  slope function and its derivatives}\label{sec:gen}
We shall measure the performance of the estimator $\widehat{\beta}_\nu$ defined in \eqref{bm:def:est:reg} by the $\cW_\nu$-risk, that is 
$\Ex\normV{\widehat{\beta}_\nu-\beta}_\nu^2$, provided that $\beta\in \cW_p$ for some $p\geq\nu\geq 0$. For an integer $\nu$ the  Sobolev norm $\normV{g}_\nu $ is equivalent to $\normV{g}+\normV{g^{(\nu)}}$, where $g^{(\nu)}$ denotes the $\nu$-th derivative of $g$ in a weak sense. Consequently, the $W_\nu$-risk reflects the performance of $\widehat{\beta}_\nu$ and $\widehat{\beta}_\nu^{(\nu)}$ as estimator of $\beta$ and $\beta^{(\nu)}$ respectively. 

The $\cW_{\nu}$-risk  is essentially determined by the deviation of the estimators of $g_k$ and $\lambda_k$, $k\in\Z$,  and  by the regularization error due to the spectral cut-off. In fact, if%
\begin{equation}\label{bm:def:est:bias}\widetilde{\beta}^{\alpha}_{\nu}:=\sum_{k\in\Z} \beta_{k} \cdot\1\{\widehat{\lambda}_{k}/  \gamma_{k}^{\nu}\geq \alpha\}\cdot \phi_{k}\text{ with }{\beta}_k:=\skalarV{{\beta},\phi_k}, \;k\in\Z,\end{equation} then by assuming $\beta\in \cW_{p}$ for some $p\geq
\nu\geq 0$ we  bound the $\cW_{\nu}$-risk of $\widehat{\beta}_\nu$ by%
\begin{equation}\label{res:dec}
 \Ex\|\widehat{\beta}_{\nu}-\beta\|_{\nu}^2\leq 2\{\Ex\normV{\widehat{\beta}_{\nu}-
 \widetilde{\beta}^{\alpha}_{\nu}}^2_{\nu}+\Ex\normV{\widetilde{\beta}^{\alpha}_{\nu}- \beta}^2_{\nu}\}.
\end{equation}
Under the moment condition $(X,W)\in\cF_{\eta,\tau}^8$ defined in \eqref{bm:def:F} and $U\in\cE^4_\eta$ we show in the proof of the next proposition that   $\Ex\normV{\widehat{\beta}_{\nu}-
 \widetilde{\beta}^{\alpha}_{\nu}}^2_{\nu}$ is  bounded up to a universal  constant by $(\alpha^2\cdot n)^{-1}\cdot\{\sigma^2 + \normV{\beta}^2 \cdot\Ex \normV{X}^2\}\cdot\Ex\normV{W}^2\cdot \eta$ and  that the 
 regularization error satisfies $\Ex\normV{\widetilde{\beta}^{\alpha}_{\nu}- \beta}^2_{\nu}=o(1)$ provided  
$\alpha=o(1)$ and $(\alpha\cdot n)^{-1}=o(1)$ as $n\to\infty$. The next assertion summarizes the minimal conditions to ensure consistency of the proposed estimator. 
 \begin{prop}[Consistency]\label{res:prop:cons}
Let $\beta\in \cW_p$, $p\geq 0$.  Consider for $0\leq \nu\leq p$ the estimator $\widehat{\beta}_\nu$ given in
(\ref{bm:def:est:reg}) with threshold  satisfying $\alpha=o(1)$ and
$(\alpha^2\cdot n)^{-1}=o(1)$ as  $n\to\infty$. If in addition $(X,W)\in\cF_\eta^8$ and  $U\in\cE^4_\eta$, then we have
$\Ex\|\widehat{\beta}_\nu-\beta\|_{\nu}^2=o(1)$ as  $n\to\infty$. 
 \end{prop}
\begin{rem}The last result covers the case $0=\nu=p$, i.e., the estimator of $\beta$ is consistent without an additional smoothness assumption on $\beta$. However,   $\widehat{\beta}_1'$ is a  consistent estimator of  $\beta'$, only if  $\beta$ is differentiable, i.e., $\beta\in\cW_p$, $p\geq1$.  \hfill$\square$\end{rem}
%%%%%%%%%%%%%%%%%%%%%%
\paragraph{Link condition.}In order to obtain a rate of convergence of 
the regularization error  and hence the
$\cW_{\nu}$-risk  we   link the smoothness condition on $\beta$, i.e., the Sobolev space  $\cW_p$, and the values $\{\lambda_k\}$ associated to the cross-covariance function of $(X,W)$.  In fact, the obtainable rate of convergence is essentially determined by the decay of $(\lambda_k)_{k\in\Z}$ as $|k|\to\infty$, which we first allow to a have  a  general form.   Notice that due to the finite   second moment of  $X$ 
the sequence  $(\lambda_k)_{k\in\Z}$  belongs to the  set $\ell^{\text{\tiny$+$}}_1$ of  nonnegative  summable sequences, i.e., $\sum_{k\in\Z} \lambda_k<\infty$. Thereby, $\lambda_+:=1\vee\max_{k\in\Z} \lambda_k$ is finite and  the rescaled sequence  $(\lambda_k/\lambda_+)_{k\in\Z}$ is taking only values in $(0,1]$. It is convenient to choose   an index function $\kappa:(0,1]\to
\R^+$ (c.f. \cite{NairPereverzevTautenhahn05}), which we always  assume here to be  a continuous and strictly increasing function with $\kappa(0+)=0$.
Then, we require that the sequence $(\lambda_k/\lambda_+)_{k\in\Z}$  is an element of  the  subset $\cS_{\kappa,d}$ of  $\ell_1^{\text{\tiny$+$}}$ defined for $d\geq 1$ by
\begin{equation}\label{res:link:gen}
\cS_{\kappa,d}:=\Bigl\{ (\lambda_k)\in\ell_1^{\text{\tiny$+$}}:\kappa\Bigl(\frac{\lambda_k}{d\,\gamma_k^{\nu}\,{\lambda}_+}\Bigr) \leq \gamma_k^{\nu-p} \leq \kappa\Bigl(\frac{d\,\lambda_k}{\gamma_k^{\nu}\,{\lambda}_+}\Bigr),\;k\in\Z\Bigr\},
\end{equation}
where  the weights $\{\gamma_k\}$ are  given in \eqref{bm:def:gamma}.  First we  consider this general class 
of values $\{\lambda_k\}$. However, we illustrate condition \eqref{res:link:gen} in Section \ref{sec:infin}  by assuming a \lq\lq regular decay\rq\rq.

The lower bound as well as the upper bound of the $\cW_{\nu}$-risk derived below involve  additional  conditions on the moments of the random function $(X,W)$, which are formalized by using the set $\cF^m_{\eta,\tau}$ defined in \eqref{bm:def:F}.   We suppose in what follows that for some index function $\kappa(\cdot)$  the random function  $(X,W)$  belongs to the subset $\cF_{\kappa}^{m}$ of $\cF^m_{\eta,\tau}$  given by%
 \begin{multline}\label{res:regressor:gen}
\cF_{\kappa}^m:=\Bigl\{ (X,W)\in\cF^m_{\eta,\tau} \text{ with  associated values }(\lambda_k)_{k\in\Z}\in  \cS_{\kappa,d}\\\text{and such that  $\Ex\normV{X}^2\leq \Lambda$, $\Ex\normV{W}^2\leq \Lambda$}\Bigr\}
\end{multline}
for some constants $d,\eta, \tau,\Lambda\geq 1$ and $m\in \N$.
%%%%%%%%%%%%%%%%%%%%%%
\paragraph{The lower bound.}
It is well-known  that in general the hardest one-dimensional subproblem does not capture the full difficulty in estimating the solution of an inverse problem even in case of a known operator (for details see e.g. the proof in \cite{MairRuymgaart96}).  In other words, there does not exist two sequences of slope functions $\beta_{1,n},\beta_{2,n}\in \cW_p^\rho$, which are statistically not consistently distinguishable and which satisfy $\normV{\beta_{1,n}-\beta_{2,n}}^2_\nu\geq C \psi_n$, where $\psi_n$ is the optimal rate of convergence. Therefore we need to consider subsets of $\cW_p^\rho$ with  growing number of elements in order to get the optimal lower bound. More specific, we obtain the following lower bound by applying Assouad's cube technique (see e.g. \cite{KorostelevTsybakov1993} or \cite{ChenReiss2011}).
\begin{theo}\label{res:lower}Assume an $n$-sample of $(Y,X,W)$ satisfying \eqref{model} and \eqref{model:W} with $\sigma>0$.  
Consider $\cW_p^\rho,$ $p,\rho>0,$  as set of slope functions,  $U\in\cE^{l}_\eta$, $ l\in\N$, as  set of error terms  and $\cF^m_\kappa$, $m\in\N$, as class of regressors defined in  \eqref{res:regressor:gen} for an arbitrary  index function $\kappa$, constants $d,\eta,\tau,\Lambda\geq 1$ and $0\leq \nu< p$. Denote by $\varphi$  the inverse function of $\kappa$. Let $k^*:=k^*(n)\in\N$ and $\delta^*:=\delta^*(n)\in(0,1]$ for some $\triangle\geq 1$  be chosen such that
\begin{equation}\label{res:lower:def:md}
\triangle^{-1}\leq \sum_{|k|\leq k^*} \frac{\gamma_{k^*}^{p-\nu}}{n\cdot \varphi(\gamma_{k}^{\nu -p})} \leq \triangle\quad\text{ and }\quad \delta^*:=\varphi(\gamma_{k^*}^{\nu-p}).
\end{equation}
If we assume in addition that  $\eta$ is sufficiently large, then
\begin{equation*} \inf_{\widetilde{\beta}} \sup_{\beta \in \cW_p^\rho, (X,W)\in\cF_{\kappa}^m,U\in\cE^{l}_\eta} \left\{ \Ex\normV{\widetilde{\beta}-\beta}^2_\nu\right\}\geq \frac{1}{4}\cdot \min\Bigl\{\frac{\sigma^2}{d\,\triangle}, \frac{\rho}{2}\Bigr\} \cdot \frac{\kappa(\delta^*)}{\Lambda}.
\end{equation*}
\end{theo}
\begin{rem} The lower bound in the last result is obtained under the assumption that the class of regressors $\cF^m_{\kappa}$ and the class of error terms $\cE^{l}_\eta$ provide a certain complexity, i.e., the uniform bound $\eta$ allows the moments of $(X,W)$ to be sufficiently large. In fact, we  ensure that for certain slope functions $\beta\in\cW_p^\rho,$ the conditional distribution  of the linear prediction error  $Y-\skalarV{\beta,\widetilde{W}}$ given the optimal linear  instrument $\widetilde{W}$  can be chosen to be Gaussian. This assumption is only needed to simplify the calculation of the distance  between distributions corresponding to different slope functions.\hfill$\square$\end{rem}
%%%%%%%%%%%%%%%%%
\paragraph{The upper bound.}
In the following theorem we provide an upper bound for the estimator $\widehat{\beta}_\nu$ defined in \eqref{bm:def:est:reg} assuming an index function $\kappa$ with the additional property that 
\begin{equation}\label{upper:kappa:cond}
\text{for all }c\geq 1:\qquad\frac{\kappa(c\cdot t)}{\kappa(t)}=O(1)\quad \text{and}\quad \frac{\kappa(t)}{\kappa(t/c)}=O(1)\quad\text{as } t\to 0.
\end{equation}
The next theorem states that the rate $\kappa(\delta^*)$ of the lower bound given  in Theorem \ref{res:lower} provides also an upper bound of 
the proposed estimator $\widehat{\beta}_\nu$. We have thus proved that  the rate $\kappa(\delta^*)$ is optimal and hence the estimator $\widehat{\beta}_\nu$   is  minimax optimal. 
\begin{theo}\label{res:upper}Assume an $n$-sample of $(Y,X,W)$ satisfying \eqref{model} and \eqref{model:W} with $\sigma>0$.  
Consider $\cW_p^\rho,$ $p,\rho>0,$  as set of slope functions,  $U\in\cE^{l}_\eta$, $ l\geq 16$, as  set of error terms  and $\cF^m_\kappa$, $m\geq 32$, as class of regressors defined in  \eqref{res:regressor:gen} for  an   index function $\kappa$  satisfying \eqref{upper:kappa:cond}, some constants $d,\eta,\tau,\Lambda\geq 1$ and $0\leq \nu< p$. Let  $\widehat{\beta}_\nu$  be  the estimator defined in \eqref{bm:def:est:reg}. 
If in addition the threshold $\alpha:=\alpha(n)$ satisfies $\alpha = 8\,d\,\Lambda\,\delta^*$, where  $\delta^*\in(0,1]$ is given in \eqref{res:lower:def:md} for some $\triangle\geq 1$, then we have
\begin{gather*}
\sup_{\beta \in W_p^\rho, (X,W)\in\cF_{\kappa}^{m},U\in\cE^{l}_\eta}\Ex\normV{\widehat{\beta}_\nu- \beta}^2_\nu\leq C \, \eta \, d \, \triangle\cdot [ \sigma^2+ \rho\, \Lambda]\cdot[\triangle\,\Lambda\,\kappa(\delta^*)+1]^4\cdot\kappa(\delta^*),
\end{gather*}
where the constant $C>0$ does only depend on the index function $\kappa$ and the constants $d,\tau,\Lambda$.
\end{theo}
\begin{rem}  We would like to stress, that for integer $\nu< p$ the Theorem \ref{res:lower} and \ref{res:upper} show together that  $\kappa(\delta^*)$ is the optimal rate of convergence for the estimation of the $\nu$-th derivative $\beta^{(\nu)}$ of $\beta$. Moreover
the $\nu$-th derivative $\widehat{\beta}_\nu^{(\nu)}$ of the in \eqref{bm:def:est:reg} proposed estimator $\widehat{\beta}_\nu$ attains this optimal rate, i.e, is minimax. \hfill$\square$\end{rem}

%% file: _4theo-prop-infin.tex
\section{Optimality  in case of a \lq\lq regular decay\rq\rq}\label{sec:infin}
In  this section we consider two special cases describing a \lq\lq regular decay\rq\rq\ of the values $\{\lambda_k\}$ associated to the cross-covariance operator $T_{_{WX}}$ of the random function $(X,W)$. In the first example we suppose the values $\{\lambda_k\}$   descend polynomial, which  in  case of a linear functional model with exogenous regressor is considered e.g. in \cite{CardotFerratySarda2003} or \cite{HallHorowitz2007}. The second example concerns  values $\{\lambda_k\}$ with exponential decay.

\paragraph{The finitely smoothing case.}Assume now the values $\{\lambda_k\}$ associated to the random function $(X,W)$  have a polynomial decay, that is\footnote{We write $a_k\asymp b_k$ if there
exists a finite positive constant $c$ such that $c^{-1}a_k\leq b_k\leq c a_k$ for all $k\in\Z$.}
\begin{equation}\label{fin:def}\lambda_k\asymp |k|^{-2a}\quad\mbox{ for some }a>0.
\end{equation}
Then straightforward calculus shows  the identity $\cR(T_{_{\widetilde{W}}})=\cW_{2a}$, where $T_{_{\widetilde{W}}}$ denotes the covariance operator associated to the optimal linear instrument $\widetilde{W}$ (see the identity \eqref{bm:e1} and its discussion in Section \ref{sec:gen}). In other words, the operator $T_{_{\widetilde{W}}}$ acts like integrating   $(2a)$-times and, hence it is called finitely smoothing.  Furthermore,  it is easily seen that%
\begin{equation}\label{res:link:pol}
\forall\; 0\leq\nu<p:\quad(\lambda_k)_{k\in\Z}\in S_{\kappa,d}\quad \text{ with } \kappa(t):=t^{(p-\nu)/(a+\nu)} \text{ and some } d\geq1.
\end{equation}
In the proof of the next proposition we shown that  the condition \eqref{res:lower:def:md} implies $\delta^*\asymp n^{-2(a+\nu)/[2(p+a)+1]}$. Thereby, we have $\kappa(\delta^*)\asymp n^{-2(p-\nu)/[2(p+a)+1]}$ and hence the lower bound   in the next assertion follows from Theorem \ref{res:lower}.
\begin{prop}\label{res:lower:pol}Let the assumptions of Theorem \ref{res:lower} be satisfied with $\kappa(t)=t^{(p-\nu)/(a+\nu)}$. Then  $\inf_{\widetilde{\beta}} \sup_{\beta \in \cW_p^\rho, (X,W)\in\cF_{\kappa}^m,U\in\cE^{l}_\eta} \bigl\{ \Ex\normV{\widetilde{\beta}-\beta}^2_\nu\bigr\}\geq C\cdot n^{-2(p-\nu)/[2(p+a)+1]}$ for some $C>0.$
\end{prop}
\noindent  On the other hand, if the threshold $\alpha$ in the definition of  the estimator $\widehat{\beta}_\nu$ given in \eqref{bm:def:est:reg} is chosen such that $\alpha\asymp n^{-2(a+\nu)/[2(p+a)+1]}$. Then by applying Theorem \ref{res:upper}  the rate $n^{-2(p-\nu)/[2(p+a)+1]}$ provides up to a constant also the upper bound of the $\cW_\nu$-risk of the  estimator $\widehat{\beta}_\nu$, which is summarized in the next proposition. We have thus proved that the rate $n^{-2(p-\nu)/[2(p+a)+1]}$ is optimal and the proposed estimator $\widehat{\beta}_\nu$ is minimax optimal. Note that the index function $\kappa(t)=t^{(p-\nu)/(a+\nu)}$ satisfies the additional condition \eqref{upper:kappa:cond}.
\begin{prop}\label{res:upper:pol}Let the assumptions of Theorem \ref{res:upper} be satisfied with $\kappa(t)=t^{(p-\nu)/(a+\nu)}$.
Consider the estimator $\widehat{\beta}_\nu$ defined in \eqref{bm:def:est:reg} with threshold $\alpha= c\cdot n^{-2(a+\nu)/[2(p+a)+1]}$, $c>0$.  Then we have $\sup_{\beta \in W_p^\rho, (X,W)\in\cF_{\kappa}^{m},U\in\cE^{l}_\eta} \bigl\{ \Ex\normV{\widehat{\beta}_\nu-\beta}^2_\nu\bigr\}=O(n^{-2(p-\nu)/[2(p+a)+1]})$.\end{prop}
\begin{rem}\label{rem:res:upper:pol:1}We shall emphasize the interesting influence of the parameters $p$ and $a$ characterizing the smoothness of $\beta$ and the smoothing property of $T_{_{\widetilde{W}}}$ respectively. As we see from Proposition \ref{res:lower:pol}  and \ref{res:upper:pol}, if the value of $a$ increases the obtainable optimal rate of convergence decreases. Therefore, the parameter $a$ is often called {\it degree of ill-posedness} (c.f. \cite{Natterer84}).
On the other hand, an increasing of the value $p$ leads to a faster optimal rate. In other words, as we expect, a smoother  slope function can be faster estimated. Finally, the estimation of higher derivatives of the slope function, i.e., increasing of the value of $\nu$, is as usual only with a slower optimal rate possible.\hfill$\square$\end{rem}
\begin{rem}\label{rem:res:upper:pol:2}There is an interesting  issue
  hidden in the parametrization we have chosen. Consider classical
  indirect regression  with known operator given by the covariance
  operator $T_{_{\widetilde{W}}}$ associated to the optimal instrument
  $\widetilde{W}$, i.e., $Y=[T_{_{\widetilde{W}}}\beta](Z)+\epsilon$
  where $Z$ has a uniform distribution on $[0,1]$ and $\epsilon$ is
  white noise  (for details see e.g. \cite{MairRuymgaart96}). Then
  given a $n$-sample of $Y$ the optimal rate of convergence of the
  $\cW_\nu$-risk of any estimator of $\beta$ is of order
  $n^{-2(p-\nu)/[2(p+2a)+1]}$,  since
  $\cR(T_{_{\widetilde{W}}})=\cW_{2a}$ (c.f. \cite{MairRuymgaart96} or
  \cite{ChenReiss2011}). However, we have shown in Proposition \ref{res:lower:pol}  and \ref{res:upper:pol} that in  FLIR the  rate  $n^{-2(p-\nu)/[2(p+a)+1]}$ is optimal. Thus comparing both rates we see that in FLIR the covariance operator $T_{_{\widetilde{W}}}$ has  the {\it degree of ill-posedness} $a$ while the same operator has in  indirect regression  a {\it degree of ill-posedness} $(2a)$. In other words in FLIR we do not face the complexity of an inversion of $T_{_{\widetilde{W}}}$ but only of its square root $T_{_{\widetilde{W}}}^{1/2}$. This, roughly speaking, may be seen as a multiplication of the  stochastic equation $Y\widetilde{W}=\skalarV{\beta,\widetilde{W}}\widetilde{W}+\epsilon\widetilde{W}$ by the inverse of $T_{_{\widetilde{W}}}^{1/2}$. Notice that $T_{_{\widetilde{W}}}$ is also the covariance operator associated to the error term  $\epsilon \widetilde{W}$. Thus multiplying to the  stochastic equation  the inverse of $T_{_{\widetilde{W}}}^{1/2}$ leads, roughly speaking, to an additive white noise and hence it is then comparable  with 
an indirect regression model with operator given by $T_{_{\widetilde{W}}}^{1/2}$. However, the operator $T_{_{\widetilde{W}}}^{1/2}$ is unknown and thus it has to  be estimated from the data.\hfill$\square$\end{rem}

\paragraph{The infinitely smoothing case.}Suppose now the values $\{\lambda_k\}$ associated to the regressors $X$ and $W$  have an exponential decay, that is 
\begin{equation}\label{res:link:exp} \lambda_k\asymp \exp(-|k|^{2a})\quad\mbox{ for some }a>0.
\end{equation}
Then it is easy to check that $\cR(T_{_{\widetilde{W}}})\subset \cW_{\nu}$ for all $\nu>0$, therefore  $T_{_{\widetilde{W}}}$
is called infinitely smoothing. In fact, the transformed values $\{|\log \lambda_k|^{-1}\}$ satisfy the polynomial condition \eqref{fin:def}. Consequently, by applying the functional calculus we have $\cR(|\log(T_{_{\widetilde{W}}})|^{-1})= \cW_{2a}$. In other words,  $|\log(T_{_{\widetilde{W}}})|^{-1}$ acts like integrating   $(2a)$-times.  Moreover,  it follows that%
\begin{equation}\label{res:link:log}
\forall\; 0\leq\nu<p:\quad(\lambda_k)_{k\in\Z}\in S_{\kappa,d}\quad \text{ with } \kappa(t):=|\log t|^{-(p-\nu)/a} \text{ and some } d\geq1.
\end{equation}
 Let $\omega$ be the inverse function of 
$\omega^{-1}(t):=t\cdot \varphi(t)$, where $\varphi$ denotes the inverse function of $\kappa$. We show  in the proof of the next proposition that in an  infinitely smoothing case the condition \eqref{res:lower:def:md} implies $1/n\asymp \delta^*\, \kappa(\delta^*)$.  Then it is straightforward to see that $\delta^*\asymp 1/(n\,\omega(1/n))$ and  $\kappa(\delta^*)\asymp \omega(1/n)$. 
Furthermore,  it is shown in \cite{Mair94} that $\omega(t)=|\log t|^{-(p-\nu)/a}(1+o(1))$ as $t \to0$. Consequently, the lower bound   in the next assertion follows again from Theorem \ref{res:lower}.
\begin{prop}\label{res:lower:log}Let the assumptions of Theorem \ref{res:lower} be satisfied with $\kappa(t)=|\log t|^{-(p-\nu)/a} $.
Then we have $\inf_{\widetilde{\beta}} \sup_{\beta \in \cW_p^\rho, (X,W)\in\cF_{\kappa}^m,U\in\cE^{l}_\eta} \bigl\{ \Ex\normV{\widetilde{\beta}-\beta}^2_\nu\bigr\}\geq C\cdot (\log n)^{-(p-\nu)/a}$ for some $C>0.$
\end{prop}
\noindent The next proposition states that  the rate $(\log n)^{-(p-\nu)/a}$ of the lower bound in Proposition \ref{res:lower:log} provides up to a constant also the upper bound of the $W_\nu$-risk of the  estimator $\widehat{\beta}_\nu$. We have thus proved that  the rate $(\log n)^{-(p-\nu)/a}$ is optimal and  $\widehat{\beta}_\nu$ is minimax-optimal. 
\begin{prop}\label{res:upper:log}Let the assumptions of Theorem \ref{res:upper} be satisfied with $\kappa(t)=|\log t|^{-(p-\nu)/a}$.
Consider the estimator $\widehat{\beta}_\nu$ defined in \eqref{bm:def:est:reg} with threshold $\alpha= c\cdot n^{-1/4}$, $c>0$. 
Then we have $\sup_{\beta \in W_p^\rho, (X,W)\in\cF_{\kappa}^{m},U\in\cE^{l}_\eta} \bigl\{ \Ex\normV{\widehat{\beta}_\nu-\beta}^2_\nu\bigr\}=O((\log n)^{-(p-\nu)/a})$.
\end{prop}
\begin{rem}It seems rather surprising that in opposite to Proposition \ref{res:upper:pol} in the last assertion the threshold  $\alpha$ does not depend on the values of $p$, $\nu$ or $a$. This, however, is due to the fact that  for $\alpha=cn^{-1/4}$,  $c>0$, the $W_\nu$-risk of $\widehat{\beta}_\nu$ is of order $O(n^{-1/2}+ |\log n^{-1/4}|^{-(p-\nu)/a})=O((\log n)^{-(p-\nu)/a})$.   Note, that the parameter $a$ specifying in condition \eqref{res:link:exp} the decay of the values $\{\lambda_k\}$ describes also in this situation the {\it degree of ill-posedness}. Finally, a comparison with an indirect regression model as in Remark \ref{rem:res:upper:pol:2} leads to the same findings.\hfill$\square$ 
\end{rem}

%% file: _5concl.tex
\section{Conclusion and perspectives}\label{sec:concl}
Assuming joint second order stationarity of the regressor $X$ and  the instrument $W$ we derive in this paper the  minimax optimal rate of convergence of an estimator of the slope function $\beta$ and its derivatives provided the covariance operator associated to the optimal linear instrument $\widetilde{W}$
is in  some general form smoothing. This results in its generality cover in particular the  case of  finitely   or infinitely smoothing covariance operators. It is worth pointing out  that for establishing the  lower bound it is not necessary to assume that the regressor and the instrument are jointly  second  order stationary. Moreover, the lower bound is derived by assuming a certain complexity of the class  of distributions  of $X$ and $W$, in particular,  it contains a Gaussian model. Therefore, we claim that replacing the 
optimal linear by the optimal instrument, i.e., the conditional expectation of $X$ given $W$, will 
not improve the optimal rate of convergence. Indeed, in a Gaussian model both instruments, if they exist, coincide.   

Many ideas in this paper can be  adapted to the general case without the assumption of joint second order stationarity of the regressor $X$ and  the instrument $W$. However, the estimation procedure itself may be different, since a projection onto the in general unknown eigenfunctions of the covariance operator of the optimal linear instrument is not possible. This is subject of ongoing research.

Once this will be established, the  open problem of how to choose the threshold $\alpha$ adaptively from the data will  remain in case not knowing the true smoothness of the slope function or not knowing the true link between the covariance operator of the optimal linear instrument and the Sobolev spaces.

%% file: _a1proofs-consistence.tex
We begin by defining and recalling notations to be used in the proofs:
\begin{multline}\label{app:l:upp:def}
X_{ik}:= \skalarV{X_i,\phi_k}, \quad W_{ik}:= \skalarV{W_i,\phi_k},\quad
T_{n,k}:=\frac{1}{n}\sum_{i=1}^n (Y_iW_{ik} - \frac{\widehat{c}_k}{\widehat{w}_k}\beta_k |{W_{ik}}|^2),\\
c_k=\Ex [\overline{X}_{ik}W_{ik}],\quad w_k=\Var(W_{ik}),\;\quad x_k=\Var(X_{ik}),\;\text{ and }\; \lambda_k=c_k^2/w_k.\end{multline} 
We shall prove in the end of this section four technical Lemma (\ref{app:l1:upp} - \ref{app:lem1:lower}) which are used in the following proofs.
\paragraph{Proof of the consistency.}\hfill\\[1ex]
\noindent\textcolor{darkred}{\sc Proof of Proposition \ref{res:prop:cons}.}
The proof is based on the  decomposition \eqref{res:dec}. We show  below  for some universal constant $C>0$ the following  bound
\begin{equation}\label{pr:res:prop:cons:e1}
\Ex\normV{\widehat{\beta}_\nu- \widetilde{\beta}^\alpha_{\nu}}^2_\nu\leq  \frac{C\,\eta}{\alpha^2 \, n}\cdot \Ex\normV{W}^{2 } \cdot\{\sigma^2+\normV{\beta}^2\Ex\normV{X}^{2}\},
\end{equation}
while in case of $\normV{\beta}_\nu<\infty$ 
we conclude from Lebesgue's dominated convergence theorem%
\begin{equation}\label{pr:res:prop:cons:e2}
\Ex\normV{\widetilde{\beta}^\alpha_{\nu}-\beta}^2_\nu=o(1)\mbox{ provided }\alpha=o(1) \mbox{ and } (\alpha\,n)^{-1}=o(1) \mbox{ as } n\to\infty.
\end{equation}
Consequently, the
conditions on  $\alpha$  ensure the convergence to zero of  the two  terms on the right hand
side in \eqref{res:dec} as $n$ tends to $\infty$, which gives the result.

Proof of (\ref{pr:res:prop:cons:e1}). By  making use of the notations given in \eqref{app:l:upp:def} it follows that
\begin{equation*}\Ex\normV{\widehat{\beta}_\nu- \widetilde{\beta}^\alpha_{\nu}}^2_{\nu}\leq\frac{1}{\alpha} \sum_{k\in\Z}  \Ex \frac{|\widehat{g}_k-\beta_k\widehat{\lambda}_k|^2}{\widehat{\lambda}_k}\1\{\widehat{\lambda}_k\geq \alpha\,\gamma_k^\nu\}\leq \frac{1}{\alpha^2} \sum_{k\in\Z}  \Ex |T_{n,k}|^2
\end{equation*}
 and hence by using \eqref{app:l1:upp:e1}  in Lemma \ref{app:l1:upp}  we obtain  (\ref{pr:res:prop:cons:e1}).

Proof of (\ref{pr:res:prop:cons:e2}). If $\beta\in W_p$, $p\geq \nu\geq 0$, then by  making use of the relation
\begin{equation*}\Ex\normV{\widetilde{\beta}^{\alpha}_{\nu}- \beta}^2_{\nu}= \sum_{k\in\Z}  \beta_k^2\cdot \gamma_k^\nu \cdot\Ex\1\{ \widehat{\lambda}_k /\gamma_k^\nu<\alpha\}\leq \sum_{k\in\Z}  \beta_k^2\cdot \gamma_k^\nu=\normV{\beta}_\nu^2\leq \normV{\beta}_p^2<\infty
\end{equation*}
the result follows from Lebesgue's dominated convergence theorem  since for each $k\in\Z$ we claim $\Ex\1\{ \widehat{\lambda}_k /\gamma_k^\nu<\alpha\}=o(1)$ provided $\alpha=o(1)$ and $(\alpha\,n)^{-1}=o(1)$  as  $n\to\infty$. Indeed, there exists $\alpha_k>0$ such that for all $\alpha\leq \alpha_k$ it holds $\lambda_k\geq 4 \tau \alpha \gamma_k^{\nu}$ and hence by using \eqref{app:l3:upp:e2}  in Lemma \ref{app:l3:upp} we bound $\Ex\1\{ \widehat{\lambda}_k /\gamma_k^\nu<\alpha\}$ up to a constant by $(\alpha\,n)^{-1} \Ex\normV{X}^2\{1+ (\alpha\,n)^{-1}\Ex\normV{W}^2\}$. Thereby, the
conditions on  $\alpha$ imply   (\ref{pr:res:prop:cons:e2}) which completes the proof. \hfill$\square$%\\
% %%%%%%%%%%%%%%%%%%%%%%%%%%%%%%%%%%%%%%%%%%

%% file: _a2proofs-lower.tex
\paragraph{Proof of the lower bound.}\hfill\\[1ex]
\noindent\textcolor{darkred}{\sc Proof of Theorem \ref{res:lower}.}  Assuming $\eta$ to be sufficiently large we can pick an 
i.i.d. sample  $(X_i,W_i)\in \cF_{\eta,\tau}^m$, $i=1,\dotsc,n,$  of 
 Gaussian random functions such that the associated sequence of values $(\lambda_k)_{k\in\Z}$ is an  element  of $\cS_{\kappa,d}$, i.e., $\{(X_i,W_i)\}\subset \cF_\kappa^m$. Consider independent error terms $\epsilon_i\sim \cN(0,1)$, $i=1,\cdots,n$, which are  independent of the random functions $\{(X_i,W_i)\}$. For $i=1,\dotsc,n$ let $\widetilde{W}_i$ be the optimal instrument given in \eqref{bm:oiv}, and denote 
${\widetilde{W}}_{ik}:=\skalarV{{\widetilde{W}}_{i},\phi_k}$, $k\in\Z$. Note, that $\widetilde{W}_{ik}$ is a centered random variable with variance $\lambda_k$.
Let  $\theta=(\theta_k)\in\{-1,1\}^{2k^*+1}$, where $k^*:=k^*(n)\in\N$ satisfies \eqref{res:lower:def:md} for some $\triangle\geq 1$. Define a  $(2k^*+1)$-vector of coefficients $(b_j)$ such that $(b_j^2)$ satisfies \eqref{app:lem1:lower:b} in Lemma \ref{app:lem1:lower}. For each $\theta$ we define a function $\beta_\theta$ which by \eqref{app:lem1:lower:e2}  in Lemma \ref{app:lem1:lower}  yields:
\begin{equation*}\beta_\theta:=\sum_{|k|\leq k^*}\theta_k b_k\phi_k\,\in \cW_p^\rho.\end{equation*}
Define  for each $\theta$ an error term $U_{\theta i}= \epsilon_i/2 + \tau_\theta\skalarV{\beta_\theta, X_i-\widetilde{W}_i}$. Then $\{U_{\theta i}\}$ are independent  centered Gaussian random variables with variance one for an appropriate chosen $\tau_\theta$, and hence  $\{U_{\theta i}\}\subset\cE_\eta^m$. Moreover, we have $\Ex[U_{\theta i}W_i(t)]=0$, for all $t\in[0,1]$ and $i=1,\dotsc,n$. Consequently, for each $\theta$ the random variables $(Y_i,X_i,W_i)$ with $Y_i:=\skalarV{\beta_{\theta},X_i}+\sigma U_{\theta i}$, $i=1,\dotsc,n,$ form a sample of the model \eqref{model}-\eqref{model:W} and we denote its joint distribution by  $P_{\theta}$. Furthermore, for $|k|\leq k^*$ and each $\theta$ we introduce $\theta^{(k)}=(\theta^{(k)}_j)\in \{-1,1\}^{2k^*+1}$ by $\theta^{(k)}_j=\theta_j$ for $k\ne j$ and $\theta^{(k)}_k=-\theta_k$. As in case of  $P_\theta$ the conditional distribution of $Y_i$ given $\widetilde{W}_i$  is Gaussian with mean  $\skalarV{\beta_\theta, \widetilde{W}_i}= \sum_{|k|\leq k^*}\theta_k b_k \overline{\widetilde{W}}_{ik}$  and variance $\sigma^2_\theta\geq \sigma^2/4$ with    $\sigma^2_\theta=\sigma^2_{\theta^{(k)}}$ it is easily seen that   the log-likelihood of $P_{\theta^{(k)}}$ w.r.t. $P_{\theta}$ is given by 
\begin{equation*}
\log\Bigl(\frac{dP_{\theta^{(k)}}}{dP_{\theta}}\Bigr)=-\frac{1}{\sigma^2_\theta}\theta_k \overline{b}_{k} \sum_{i=1}^n 
U_{\theta i}\widetilde{W}_{ik} -\frac{1}{\sigma^2_\theta}\theta_k b_{k} \sum_{i=1}^n 
\overline{U}_{\theta i}\overline{\widetilde{W}}_{ik}+\frac{2}{\sigma^2_\theta}b_k^2 \sum_{i=1}^n \widetilde{W}^2_{ik}.
\end{equation*}
Its expectation satisfies 
%\begin{equation*}
$\Ex_{P_{\theta}}[\log(dP_{\theta^{(k)}}/dP_{\theta})]\geq - n\cdot b^2_k\cdot \lambda_k / (2\sigma^2)$ since $\lambda_k=\Var \widetilde{W}_{ik}$ and $\sigma^2_\theta\geq \sigma^2/4$. 
%\end{equation*}
In terms of  Kullback-Leibler divergence this means $KL(P_{\theta^{(k)}},P_{\theta})\leq n\cdot b^2_k\cdot \lambda_k /(2\sigma^2)$. Since the
 Hellinger distance $H(P_{\theta^{(k)}},P_{\theta})$ satisfies $H^2(P_{\theta^{(k)}},P_{\theta}) \leq KL(P_{\theta^{(k)}},P_{\theta})$  it follows from \eqref{app:lem1:lower:e1} in Lemma \ref{app:lem1:lower}  that 
\begin{equation}\label{pr:lower:e3}
H^2(P_{\theta^{(k)}},P_{\theta}) \leq \frac{n }{2\sigma^2}\cdot b^2_k\cdot \lambda_k\leq 1,\quad |k|\leq k^*.
\end{equation} 
Consider  the  Hellinger affinity $\rho(P_{\theta^{(jk}},P_{\theta})= \int \sqrt{dP_{\theta^{(k)}}dP_{\theta}}$, then we obtain for any estimator $\widetilde\beta$  that
\begin{align}\nonumber
\rho(P_{\theta^{(k)}},P_{\theta})&\leq \int \frac{|\skalarV{\widetilde{\beta}-\beta_{\theta^{(k)}},\phi_k}|}{|\skalarV{\beta_{\theta}-\beta_{\theta^{(k)}},\phi_k}|} \sqrt{dP_{\theta^{(k)}}dP_{\theta}} + \int \frac{|\skalarV{\widetilde{\beta}-\beta_{\theta},\phi_k}|}{|\skalarV{\beta_{\theta}-\beta_{\theta^{(k)}},\phi_k}|} \sqrt{ dP_{\theta^{(k)}}dP_{\theta}}\\\label{pr:lower:e4}
&\leq \Bigl( \int  \frac{|\skalarV{\widetilde{\beta}-\beta_{\theta^{(k)}},\phi_k}|^2}{|\skalarV{\beta_{\theta}-\beta_{\theta^{(k)}},\phi_k}|^2} dP_{\theta^{(k)}}\Bigr)^{1/2} +   \Bigl( \int  \frac{|\skalarV{\widetilde{\beta}-\beta_{\theta},\phi_k}|^2}{|\skalarV{\beta_{\theta}-\beta_{\theta^{(k)}},\phi_k}|} dP_{\theta}\Bigr)^{1/2}.
\end{align}
Due to the identity $\rho^2(P_{\theta^{(k)}},P_{\theta})=1-\frac{1}{2}H^2(P_{\theta^{(k)}},P_{\theta})$   combining  \eqref{pr:lower:e3} with 
 \eqref{pr:lower:e4} yields
\begin{equation*}
\Bigl\{\Ex_{{\theta^{(k)}}}|\skalarV{\widetilde{\beta}-\beta_{\theta^{(k)}},\phi_j}|^2+ \Ex_{{\theta}}|\skalarV{\widetilde{\beta}-\beta_{\theta},\phi_k}|^2\Bigr\}\geq\frac{1}{4} b_k^2,\quad |k|\leq k^*.
 \end{equation*}
From this  we conclude for each estimator $\widetilde\beta$ that
\begin{align*}
\sup_{\beta \in \cW_p^\rho, U \in \cE_\eta^l,\atop (X,W)\in\cF_{\kappa}^m} &\Ex\normV{\widetilde\beta -\beta}_\nu^2 \geq \sup_{\theta\in \{-1,1\}^{2k^*+1}} \Ex_\theta\normV{\widetilde\beta -\beta_\theta}_\nu^2\\
&\geq \frac{1}{2^{{2k^*+1}}}\sum_{\theta\in \{-1,1\}^{2k^*+1}}\sum_{|k|\leq k^*}\gamma_k^\nu\cdot\Ex_{{\theta}}|\skalarV{\widetilde{\beta}-\beta_{\theta},\phi_k}|^2\\
&= \frac{1}{2^{{2k^*+1}}}\sum_{\theta\in \{-1,1\}^{2k^*+1}}\sum_{|k|\leq k^*}\frac{\gamma_k^\nu}{2}\cdot\Bigl\{\Ex_{{\theta}}|\skalarV{\widetilde{\beta}-\beta_{\theta},\phi_k}|^2+\Ex_{{\theta^{(k)}}}|\skalarV{\widetilde{\beta}-\beta_{\theta^{(k)}},\phi_k}|^2 \Bigr\}\\
&\geq\frac{1}{8} \sum_{|k|\leq k^*} \gamma_k^\nu b_k^2 \geq \frac{1}{4}\cdot \min\Bigl\{\frac{\sigma^2}{d\,\triangle}, \frac{\rho}{2}\Bigr\} \cdot \frac{\kappa(\delta^*)}{\Lambda},
\end{align*}
where the last inequality follows from  \eqref{app:lem1:lower:e3} in Lemma \ref{app:lem1:lower} together with $\lambda_+\leq \Lambda$ which completes the proof.\hfill$\square$%\\

%% file: _a3proofs-upper.tex
\paragraph{Proof of the upper bound.}\hfill\\[1ex]
\noindent\textcolor{darkred}{\sc Proof of Theorem \ref{res:upper}.}
The proof is based on the  decomposition \eqref{res:dec}, where  we show  below under the condition $(X,W)\in \cF_{\eta,\tau}^{32}$, $U\in \cE_\eta^{16}$ and $(\lambda_k)_{k\in\Z} \in\cS_{\kappa,d}$ for some universal constant $C>0$ and  $I:=\{k\in\Z:8\,\lambda_k> \alpha\,\gamma_k^\nu\}$ the following two bounds
\begin{multline}\label{pr:res:upper:gen:e1}
\Ex\normV{\widehat{\beta}_\nu- \widetilde{\beta}^\alpha_{\nu}}^2_\nu
\leq C\cdot \eta\cdot \Bigl\{ \frac{ 1}{\alpha\,n} + d\cdot  \sum_{k\in I}\frac{1}{n \cdot\varphi(\gamma_k^{\nu-p})}\Bigr\}\cdot\\ \cdot\{ \normV{\beta}^{2}\cdot \Ex\normV{X}^2+\sigma^{2}\}\cdot  \Bigl\{ \frac{\Ex\normV{X}^2}{\alpha\, n} +1\Bigr\}\cdot \Bigl\{ \frac{\Ex\normV{W}^2}{\alpha\, n} +1\Bigr\}^3,
\end{multline}
\begin{equation}\label{pr:res:upper:gen:e2}
\Ex\normV{\widetilde{\beta}^\alpha_{\nu}-\beta}^2_\nu\leq  \kappa(d\,4\,\tau\, \alpha)\cdot \normV{\beta}_p^2 + \frac{C\,\eta }{\alpha\, n}\cdot\Ex\normV{X}^2\cdot \normV{\beta}_p^2 \cdot\Bigl(1+\frac{\Ex\normV{W}^2}{\alpha\, n}\Bigr).
\end{equation}
Consequently, for all $\beta\in \cW_p^\rho$ and $(X,W)\in \cF^{32}_\kappa$, i.e.,  $\Ex\normV{X}^2\leq \Lambda$, $\Ex\normV{W}^2\leq \Lambda$, and hence $\lambda_+\leq \Lambda$, follows
\begin{multline*}
\Ex\normV{\widehat{\beta}_\nu- \beta}^2_\nu\leq C \cdot \eta\cdot \Bigl\{ \frac{ 1}{\alpha\,n} + d\cdot  \sum_{k\in I}\frac{1}{n\cdot \varphi(\gamma_k^{\nu-p})}+\kappa(d\,4\,\tau\, \alpha)\Bigr\}\cdot[ \sigma^2+ \rho\cdot \Lambda]\cdot \Bigl[\frac{\Lambda}{\alpha\,n}+1\Bigr]^4 .
\end{multline*}
Let $k^*:=k^*(n)\in\N$ and $\delta^*:=\delta^*(n)\in(0,1]$ be given by \eqref{res:lower:def:md} for some $\triangle\geq 1$ then the condition on  $\alpha$, that is $\alpha = 8\,d\,\Lambda\,\delta^*$, implies $I\subset\{k\in\Z:|k|\leq k^*\}$. We conclude from \eqref{res:lower:def:md} that $1/[\alpha\cdot n] \leq \triangle\cdot\kappa(\delta^*)$ and $\sum_{k\in I}1/[n\cdot\varphi(\gamma_k^{\nu-p})]\leq \triangle\cdot\kappa(\delta^*)$, hence that
\begin{multline*}
\Ex\normV{\widehat{\beta}_\nu- \beta}^2_\nu\leq C \, \eta\, d\, \triangle\cdot\frac{\kappa(d\,4\,\tau\, \alpha)}{\kappa(\alpha/(8\,d\,\Lambda))}\cdot[ \sigma^2+ \rho\, \Lambda]\cdot \Bigl[\triangle\,\Lambda\,\kappa(\delta^*)+1\Bigr]^4\cdot \kappa(\delta^*) .
\end{multline*}
Thereby, the condition \eqref{upper:kappa:cond}, that is  $\kappa(d\,4\,\tau\, \alpha)/\kappa(\alpha)=O(1)$ and $\kappa(\alpha)/\kappa(\alpha/(8\,d\,\Lambda))=O(1)$, as $\alpha$ tends to zero, implies  the result.

Proof of (\ref{pr:res:upper:gen:e1}). By using  $T_{n,k}$ introduced in \eqref{app:l:upp:def} we  obtain the identity
 \begin{equation*}
\Ex\normV{\widehat{\beta}_\nu- \widetilde{\beta}^\alpha_{\nu}}^2_\nu
=  \sum_{k\in\Z} \gamma_k^\nu \cdot \Ex\Bigl[\frac{|T_{n,k}|^2|\widehat{c}_k/\widehat{w}_k|^2\1\{\widehat{w}_k\geq\alpha\}}{
\widehat{\lambda}_k^2} \1\{\widehat{\lambda}_k\geq \alpha\, \gamma_k^\nu\}\Bigr].
 \end{equation*}
We partition the sum into two parts which we estimate separately using the bounds in Lemma \ref{app:l2:upp}. 
First by using  
$  \widehat{\lambda}_k \geq \alpha \gamma_k^\nu$ together with \eqref{app:l2:upp:e4} in Lemma \ref{app:l2:upp} we bound the sum over  $I^c:=\{k\in\Z:8\,\lambda_k\leq \alpha\,\gamma_k^\nu\}$ by 
\begin{equation*}
\sum_{k\in I^c} \Ex\Bigl[\frac{|T_{n,k}|^2|\widehat{c}_k|^2}{
\alpha^2\,\widehat{w}_k^2} \1\{\widehat{w}_k\geq\alpha\}\1\{\widehat{\lambda}_k\geq \alpha\, \gamma_k^\nu\}\Bigr]
\leq  \frac{C\,\eta \,\Ex\normV{X}^2}{n^2\,\alpha^2}\cdot\{ \normV{\beta}^{2}\cdot \Ex\normV{X}^2+\sigma^{2}\}\cdot \Bigl\{ \frac{\Ex\normV{W}^2}{\alpha\, n} +1\Bigr\}^3.
 \end{equation*}
While due to the identity $\widehat{\lambda}_k=|\widehat{c}_k|^2/\widehat{w}_k\1\{\widehat{w}_k\geq\alpha\}$
together with \eqref{app:l2:upp:e3} in Lemma \ref{app:l2:upp} the sum  over  $I:=\{k\in\Z:8\,\lambda_k> \alpha\,\gamma_k^\nu\}$ is bounded by 
\begin{multline*}
\sum_{k\in I} \gamma_k^\nu\cdot\Ex\Bigl[\frac{|T_{n,m}|^2}{\widehat{w}_k\cdot \widehat{\lambda}_k} \1\{\widehat{w}_k\geq\alpha\}\1\{\widehat{\lambda}_k\geq \alpha\, \gamma_k^\nu\}\Bigr]
\\\leq \sum_{k\in I}\frac{C\,\eta\,\gamma_k^\nu}{n\,\lambda_k}\cdot
\{\normV{\beta}^{2}\cdot \Ex\normV{X}^2+\sigma^{2}\}\cdot  \Bigl\{ \frac{\Ex\normV{X}^2}{\alpha\, n} +1\Bigr\}\cdot  \Bigl\{ \frac{\Ex\normV{W}^2}{\alpha\, n} +1\Bigr\}.
 \end{multline*}
From $\lambda_+\geq 1$ it follows that by combining the two parts of the sum we have 
\begin{multline*}
\Ex\normV{\widehat{\beta}_\nu- \widetilde{\beta}^\alpha_{\nu}}^2_\nu
\leq C\cdot \eta\cdot \Bigl\{ \frac{ 1}{\alpha\,n} + \sum_{k\in I^c}\frac{\gamma_k^\nu\,\lambda_+}{\lambda_k\, n}\Bigr\}\cdot\\ \cdot\{ \normV{\beta}^{2}\cdot \Ex\normV{X}^2+\sigma^{2}\}\cdot  \Bigl\{ \frac{\Ex\normV{X}^2}{\alpha\, n} +1\Bigr\}\cdot \Bigl\{ \frac{\Ex\normV{W}^2}{\alpha\, n} +1\Bigr\}^3.
 \end{multline*}
Now the link condition $(\lambda_k)_{k\in\Z}\in\cS_{\kappa,d}$ implies \eqref{pr:res:upper:gen:e1}.

The proof of (\ref{pr:res:upper:gen:e2}) is based on the identity
$\Ex\normV{\widetilde{\beta}^\alpha_{\nu}-\beta}^2_\nu= \sum_{k\in\Z}\gamma_k^\nu \beta_k^2 P(\widehat{\lambda}_k< \alpha\, \gamma_k^\nu),$
where we  partition the sum again into two parts which we estimate separately. First we sum over  $I:=\{k\in\Z:\lambda_k< 4\,\tau\,\alpha\,\gamma_k\}$.  Since $\lambda_+\geq 1$, the link condition $(\lambda_k)_k\in\cS_{\kappa,d}$ together with the  monotonicity of $\kappa$ shows that
\begin{equation*}
\sum_{k\in I}\gamma_k^\nu \beta_k^2 P(\widehat{\lambda}_k< \alpha\, \gamma_k^\nu)\leq
\sum_{k\in I}\gamma_k^p \beta_k^2 \kappa(d\frac{\lambda_k}{\gamma_k^\nu\lambda_+})\leq   \kappa(d\,4\,\tau\,\alpha)\cdot\sum_{k\in I}\gamma_k^p \beta_k^2.
 \end{equation*}
The sum over $I^c:=\{k\in\Z:\lambda_k\geq  4\,\tau\,\alpha\,\gamma_k\}$ we bound using  \eqref{app:l3:upp:e2}  in Lemma \ref{app:l3:upp}, that is
\begin{equation*}
\sum_{k\in I^c}\gamma_k^\nu \beta_k^2 P(\widehat{\lambda}_k< \alpha\, \gamma_k^\nu)\leq
\frac{C\,\eta}{\alpha\, n} \cdot\Ex\normV{X}^2\cdot\Bigl\{1+\frac{\Ex\normV{W}^2}{\alpha\, n}\Bigr\} \sum_{k\in I^c}\gamma_k^\nu \beta_k^2.
 \end{equation*}
Combining the two parts of the sum we obtain (\ref{pr:res:upper:gen:e2}), which completes the proof. \hfill$\square$%\\
% %%%%%%%%%%%%%%%%%%%%%%%%%%%%%%%%%%%%%%%%%%
 

%% file: _a4proofs-technical.tex
\paragraph{Technical assertions.}\hfill\\[1ex]
The following four lemma gather technical results used in the proof of Proposition \ref{res:prop:cons}, Theorem \ref{res:lower} and Theorem \ref{res:upper}. 
\begin{lem}\label{app:l1:upp} Suppose $(X,W)\in\cF_{\eta,\tau}^{4m}$ and $U\in\cE^{4m}_\eta$, $m\in\N$. Then for some constant $C>0$ only depending on $m$ we have
\begin{align}\label{app:l1:upp:e1}
\sup_{k\in\Z}\Bigl\{\frac{1}{w_k^{m}}\cdot  \Ex |T_{n,k}|^{2m}\Bigr\} &\leq C \cdot \frac{1}{n^m}\cdot\{ \normV{\beta}^{2m}\cdot (\Ex\normV{X}^2)^m+\sigma^{2m}\}\cdot \eta,\\
\label{app:l1:upp:e2}
\sup_{k\in\Z}\Ex \Bigl|\frac{\widehat{w}_k-w_k}{w_k}\Bigr|^{2m}&\leq C\cdot \frac{1}{n^m}\cdot \eta,\\
\label{app:l1:upp:e3}
\sup_{k\in\Z}\Bigl\{\frac{\lambda_k^m}{x_k^{m}}\cdot \Ex \Bigl|\frac{\widehat{c}_k-c_k}{c_k}\Bigr|^{2m}\Bigr\}&\leq C\cdot \frac{1}{n^m}\cdot\eta.
\end{align}
\end{lem}
\noindent\textcolor{darkred}{\sc Proof.} Let $\eta_{ik}:=\sum_{l\ne k }\beta_l \overline{X}_{il}$, $\zeta_{ik}:=\beta_k\{\overline{X}_{ik}-\overline{W}_{ik}c_k/w_k\}$, $\tau_{ik}:=\beta_k\overline{W}_{ik}\{c_k/w_k-\widehat{c}_k/\widehat{w}_k\}$, $i=1,\dotsc,n$ and $k\in\Z$. Then we have
\begin{equation*}T_{n,k}=\frac{1}{n}\sum_{i=1}^n\{ \eta_{ik}+\zeta_{ik}+\tau_{ik}+ \sigma U_i\}W_{ik}=:T_1+T_2+T_3+T_4,\end{equation*}
where we bound below each summand separately, that is
\begin{align}\label{app:l1:upp:e1:1}
\Ex|T_j|^{2m}&\leq C\cdot \frac{w_k^m}{n^{m}}\cdot  \normV{\beta}^{2m} \cdot  (\Ex\normV{X}^2)^m\cdot \eta,\quad j=1,2,3,\\\label{app:l1:upp:e1:2}
\Ex|T_4|^{2m}&\leq C\cdot \frac{w_k^m}{n^{m}}\cdot\sigma^{2m} \cdot \eta
\end{align}
for some $C>0$ only depending on $m$. Consequently, the  inequality \eqref{app:l1:upp:e1} follows from \eqref{app:l1:upp:e1:1} and \eqref{app:l1:upp:e1:2}. Consider $T_1$. For each $k\in \Z$ the random variables $(\eta_{ik}\cdot W_{ik})$,  $i=1,\dots,n,$ are independent and identically distributed with mean zero. From Theorem 2.10 in \cite{Petrov1995} we conclude $\Ex |T_{1}|^{2m} \leq C n^{-m}\Ex |\eta_{1k} W_{1k}|^{2m}$ for some  constant $C>0$ only depending on $m$. Then  we claim that \eqref{app:l1:upp:e1:1} follows in case of  $T_1$ from  the Cauchy-Schwarz inequality together with   $(X,W)\in\cF_{\eta,\tau}^{4m}$, i.e.,  $\sup_{i,k} \Ex|X_{ik}/\sqrt{x_k}|^{4m}\leq \eta$ and  $\sup_{i,k} \Ex|W_{ik}/\sqrt{w_k}|^{4m}\leq \eta$. Indeed, we have
$$\Ex |\eta_{1k}W_{1k}|^{2m}
\leq (\sum_{j\ne k}\beta_j^2)^m\sum_{j_1\neq k}\dots\sum_{j_m\ne k} \Ex |W_{ik}|^{2m}\prod_{l=1}^m |X_{1j_l}|^2\leq \normV{\beta}^{2m}\cdot w_k^m \cdot  (\sum_{j\neq k}x_{j})^m\cdot \eta.$$
Consider $T_2$. \eqref{app:l1:upp:e1:1} follows  in analogy to the case of $T_1$ since $\{\zeta_{ik}\cdot W_{ik} \}$ are independent and identically distributed with mean zero respectively, and 
$\Ex |\zeta_{1k}\cdot W_{1k}|^{2m}\leq C\cdot \beta_k^{2m}\{ x_k^m \cdot w_k^m \cdot \eta + |c_k|^{2m} \eta\} \leq C\cdot
\normV{\beta}^{2m} \cdot (\Ex\normV{X}^2)^m\cdot w_k^m\cdot \eta,$ where  $C>0$ does only depend on $m$. Consider $T_3$. We have $\Ex |T_{3}|^{2m}\leq C \beta_k^{2m} \{ |c_k|^{2m}\Ex|\widehat{w}_k/w_k-1|^{2m} + \Ex|\widehat{c}_k-c_k|^{2m}\}$ 
for some $C>0$ only depending on $m$, by the identity $T_3= \beta_k\{\widehat{w}_kc_k/w_k - \widehat{c}_k\}$. Therefore \eqref{app:l1:upp:e1:1} in case of  $T_3$ follows from \eqref{app:l1:upp:e2} and \eqref{app:l1:upp:e3}. Consider $T_4$. \eqref{app:l1:upp:e1:2} follows in analogy to the case of $T_1$, because  $\{\sigma U_{i}\cdot W_{ik} \}$ are independent and identically distributed with mean zero respectively, and 
$\Ex |\sigma \cdot U_1\cdot W_{1k}|^{2m}\leq C\cdot \sigma^{2m} \cdot w_k^m \cdot \eta,$ where  $C>0$ does only depend on $m$.

Proof of (\ref{app:l1:upp:e2}) and (\ref{app:l1:upp:e3}).  Since $\{(|W_{ik}|^2/w_k-1)\}$  and $\{(\overline{X}_{ik}W_{ik}-c_k)\}$ are independent and identically distributed with mean zero respectively, where 
$\Ex| |W_{1k}|^2/w_k|^{2m}\leq \eta$ and $\Ex| \overline{X}_{1k}W_{1k}-c_k|^{2m}\leq C\cdot x_k^m \cdot w_k^m \cdot \eta $,  for some $C>0$ only depending on $m$, the result follows by applying Theorem 2.10 in \cite{Petrov1995}, which proves the lemma.\hfill$\square$\\
%%%%%%%%%%%%%%%%%%%%%%%%%%%%%%%%%%%%%%%%%%%
%%%%%%%%%%%%%%%%%%%%%%%%%%%%%%%%%%%%%%%%%%%
%%%%%%%%%%%%%%%%%%%%%%%%%%%%%%%%%%%%%%%%%%%
\begin{lem}\label{app:l3:upp} Let $(X,W)\in\cF_{\eta,\tau}^{4m}$, $m\in\N$, then for all $0<d<1$ and some  constant $C>0$ only depending on $m$ we have
\begin{equation}\label{app:l3:upp:e1} \sup_{k\in\Z}P(\widehat{w}_k/w_k<d )\leq C \cdot \frac{1}{(1-d)^{2m}}\cdot \frac{1}{n^m}\cdot \eta.
\end{equation}
Suppose $(X,W)\in\cF_{\eta,\tau}^{8}$ and let $I_1:=\{k\in\Z: \lambda_k\geq 4\,\tau \, \alpha\,\gamma_k^\nu\}$. Then for some universal constant $C>0$ we have 
\begin{equation}\label{app:l3:upp:e2}\sup_{k\in I_1}P(\widehat{\lambda}_k<\alpha\, \gamma_k^\nu )\leq \frac{C\,\eta}{\alpha\, n} \cdot\Ex\normV{X}^2\cdot\Bigl\{1+\frac{\Ex\normV{W}^2}{\alpha\, n}\Bigr\} .
\end{equation}
While if $(X,W)\in\cF_{\eta,\tau}^{16}$ and  $I_2:=\{k\in\Z: 8\,\lambda_k\leq \alpha\,\gamma_k^\nu\}$. Then
\begin{equation}\label{app:l3:upp:e3}\sup_{k\in I_2}\Bigl\{\frac{\lambda_k^2}{x_k^2}\cdot P(\widehat{\lambda}_k\geq\alpha\, \gamma_k^\nu )\Bigr\}\leq \frac{C\,\eta}{n^2}\cdot \Bigl\{1+\frac{\Ex\normV{W}^2}{\alpha\, n}\Bigr\}^2.
\end{equation}
\end{lem}
\noindent\textcolor{darkred}{\sc Proof.} Since  $P(\widehat{w}_k/w_k< d)\leq P(|\widehat{w}_k/w_k-1|\geq 1-d)$ by  applying Markov's inequality the estimate \eqref{app:l3:upp:e1} follows  from \eqref{app:l1:upp:e2} in Lemma \ref{app:l1:upp}.

The proof of \eqref{app:l3:upp:e2} is based on the  elementary inequality
$$P(\widehat{\lambda}_k<\alpha\,\gamma_k^\nu) \leq P(\widehat{w}_k<\alpha) + P(\widehat{\lambda}_k<\alpha\,\gamma_k^\nu \wedge \widehat{w}_k\geq\alpha),$$ where we show below for some universal constant $C>0$ the following two bounds
\begin{gather}\label{app:l3:upp:e2:1}
\sup_{k\in I_1}P(\widehat{w}_k<\alpha) \leq C \cdot \frac{1}{n}\cdot \eta,\\\label{app:l3:upp:e2:2}
\sup_{k\in I_1}P(\widehat{\lambda}_k<\alpha\,\gamma_k^\nu \wedge \widehat{w}_k\geq\alpha)\leq \frac{C\eta}{\alpha n} \cdot\Ex\normV{X}^2\cdot\Bigl\{1+\frac{\Ex\normV{W}^2}{\alpha\, n} \Bigr\} 
\end{gather}
which imply together \eqref{app:l3:upp:e2}.  Since $\alpha/w_k \leq \lambda_k /w_k \cdot \alpha \,\gamma_k^\nu/\lambda_k\leq \tau \alpha \,\gamma_k^\nu/\lambda_k\leq 1/4$ holds true for all $k\in I_1$, the estimate \eqref{app:l3:upp:e2:1} follows from \eqref{app:l3:upp:e1}. The proof of \eqref{app:l3:upp:e2:2} is based on 
\begin{equation}\label{app:l3:upp:e2:3} 
1-2\frac{\widehat{\lambda}_k}{\lambda_k}\leq \Bigl\{4|\widehat{c}_k/c_k-1|^2+1\Bigr\}\cdot\Bigl\{ \frac{|\widehat{w}_k/w_k-1|^2}{\widehat{w}_k} +{|\widehat{w}_k/w_k-1|}\Bigr\}+ 4|\widehat{c}_k/c_k-1|^2\end{equation}
which implies for all $k\in I_1$ that 
\begin{multline*}
P(\widehat{\lambda}_k<\alpha\,\gamma_k^\nu \wedge \widehat{w}_k\geq\alpha)\leq P\Bigl(1/4\leq 4|\widehat{c}_k/c_k-1|^2\Bigr)\\+P\Bigl(
1/4\leq \Bigl\{4|\widehat{c}_k/c_k-1|^2+1\Bigr\}\cdot\Bigl\{ \frac{|\widehat{w}_k/w_k-1|^2}{\alpha} +{|\widehat{w}_k/w_k-1|}\Bigr\}\Bigr).\end{multline*}
Therefore, by applying  Markov's inequality   together with (\ref{app:l1:upp:e2}) and (\ref{app:l1:upp:e3}) in Lemma \ref{app:l1:upp} we obtain the estimate \eqref{app:l3:upp:e2:2}.  

The proof of \eqref{app:l3:upp:e3} is based on the  decomposition
\begin{align*}
P(\widehat{\lambda}_k\geq\alpha\, \gamma_k^\nu )&\leq P\Bigl(|\widehat{c}_k/c_k-1|^2\frac{\1\{\widehat{w}_k\geq\alpha\}}{\widehat{w}_k/w_k}\geq \frac{\alpha\gamma_k^\nu}{4\lambda_k}\Bigr) + P\Bigl(\frac{\1\{\widehat{w}_k\geq\alpha\}}{\widehat{w}_k/w_k} \geq \frac{\alpha\gamma_k^\nu}{4\lambda_k}\Bigr)
\end{align*}
which implies for all $k\in I_2$ together with Markov's inequality that 
\begin{align*}
P(\widehat{\lambda}_k\geq\alpha\, \gamma_k^\nu )
&\leq \frac{1}{4} \Ex\Bigl[|\widehat{c}_k/c_k-1|^4\frac{\1\{\widehat{w}_k\geq\alpha\}}{|\widehat{w}_k/w_k|^2}\Bigr] + P\Bigl({\widehat{w}_k/w_k} \leq 1/2\Bigr).
\end{align*}
Therefore, by using   
$ 1\leq 2^3\{|\widehat{w}_k/w_k-1|^4+ |\widehat{w}_k/w_k|^2|\widehat{w}_k/w_k-1|^2+|\widehat{w}_k/w_k|^2\}$ 
we obtain
\begin{align*}
P(\widehat{\lambda}_k\geq\alpha\, \gamma_k^\nu )
&\leq 4\Ex \Bigl[|\widehat{c}_k/c_k-1|^4 \{\frac{|\widehat{w}_k/w_k-1|^4}{\alpha^2/w_k^2} +|\widehat{w}_k/w_k-1|^2 +1\}\Bigr] + P\Bigl({\widehat{w}_k/w_k} \leq 1/2\Bigr).
\end{align*}
Now (\ref{app:l1:upp:e2}) and (\ref{app:l1:upp:e3}) in Lemma \ref{app:l1:upp} and \eqref{app:l3:upp:e1} gives  \eqref{app:l3:upp:e3}, which completes the proof. \hfill$\square$\\
%%%%%%%%%%%%%%%%%%%%%%%%%%%%%%%%%%%%%%%%%%%
%%%%%%%%%%%%%%%%%%%%%%%%%%%%%%%%%%%%%%%%%%%
%%%%%%%%%%%%%%%%%%%%%%%%%%%%%%%%%%%%%%%%%%%
\begin{lem}\label{app:l2:upp} Suppose $(X,W)\in\cF_{\eta,\tau}^{32}$ and $U\in\cE^{16}_\eta$. Let $I:=\{k\in \Z:8\,\lambda_k\leq \alpha\,\gamma^\nu\}$, then for some universal constant $C>0$  we have
\begin{multline}\label{app:l2:upp:e1}
\sup_{k\in\Z}\,\Bigl\{ \lambda_k\cdot \Ex \Bigl[\frac{|T_{n,k}|^{2}{|\widehat{c}_k/c_k-1|^2}\1\{\widehat{w}_k\geq \alpha\}}{\widehat{w}_k}\Bigr]\Bigr\} \\\leq  \frac{C\,\eta}{n}\cdot
\{\normV{\beta}^{2}\cdot \Ex\normV{X}^2+\sigma^{2}\}\cdot \frac{\Ex\normV{X}^2}{n}\cdot 
\Bigl\{ \frac{\Ex\normV{W}^2}{\alpha\, n} +1\Bigr\},\end{multline}
\begin{multline}\label{app:l2:upp:e2}
\sup_{k\in\Z}\,\Bigl\{ w_k^2\cdot \Ex \Bigl[\frac{|T_{n,k}|^{4}\1\{\widehat{w}_k\geq \alpha\}}{\widehat{w}_k^4}\Bigr]\Bigr\}\\ \leq  \frac{C\,\eta}{n^2}\cdot\{ \normV{\beta}^{2}\cdot \Ex\normV{X}^2+\sigma^{2}\}^2\cdot \Bigl\{ \frac{\Ex\normV{W}^2}{\alpha\, n} +1\Bigr\}^4,
\end{multline}
\begin{multline}\label{app:l2:upp:e3}
\sup_{k\in\Z}\Bigl\{\lambda_k\cdot\Ex\Bigl[\frac{|T_{n,k}|^2}{\widehat{w}_k\cdot \widehat{\lambda}_k} \1\{\widehat{w}_k\geq\alpha\}\1\{\widehat{\lambda}_k\geq \alpha\, \gamma_k^\nu\}\Bigr]\Bigr\}\\\leq 
\frac{C\,\eta}{n}\cdot\{ \normV{\beta}^{2}\cdot \Ex\normV{X}^2+\sigma^{2}\}\cdot \Bigl\{ \frac{\Ex\normV{X}^2}{\alpha\, n} +1\Bigr\}\cdot \Bigl\{ \frac{\Ex\normV{W}^2}{\alpha\, n} +1\Bigr\},
 \end{multline}
\begin{multline}\label{app:l2:upp:e4}
\sup_{k\in I}\Bigl\{\frac{1}{x_k}\cdot \Ex \Bigl[\frac{|T_{n,k}|^2\,|\widehat{c}_k|^2}{\widehat{w}_k^2}\1\{\widehat{w}_k\geq \alpha\}\1\{\widehat{\lambda}_k\geq \alpha\,\gamma_k^\nu\}\Bigr]\Bigr\} \\\leq  \frac{C\,\eta}{n^2}\cdot\{ \normV{\beta}^{2}\cdot \Ex\normV{X}^2+\sigma^{2}\}\cdot \Bigl\{ \frac{\Ex\normV{W}^2}{\alpha\, n} +1\Bigr\}^3.
\end{multline}
\end{lem}
\noindent\textcolor{darkred}{\sc Proof.} Consider the elementary inequality \begin{equation}\label{app:l2:upp:e1:1} 1\leq 2\Bigl\{|\widehat{w}_k/w_k-1|^2+ |\widehat{w}_k/w_k||\widehat{w}_k/w_k-1|+|\widehat{w}_k/w_k|\Bigr\}.\end{equation}By applying the Cauchy-Schwarz inequality it follows that
\begin{multline*}
\Ex \Bigl[\frac{|T_{n,k}|^{2}|\widehat{c}_k/c_k-1|^2\1\{\widehat{w}_k\geq \alpha\}}{\widehat{w}_k}\Bigr] \\
\leq 2 \Bigl(\Ex|T_{n,k}|^4\Bigl)^{1/2}\Bigl(
\Ex  |\widehat{c}_k/c_k-1|^8\Bigr)^{1/4}\Bigl\{\frac{(\Ex|\widehat{w}_k/w_k-1|^8)^{1/2}}{\alpha^2} +\frac{(\Ex|\widehat{w}_k/w_k-1|^4)^{1/2}}{w_k^2}+\frac{1}{w_k^2}\Bigr\}^{1/2}
\end{multline*}
and, hence  \eqref{app:l1:upp:e1}, \eqref{app:l1:upp:e2} and \eqref{app:l1:upp:e3} in Lemma \ref{app:l1:upp} imply  \eqref{app:l2:upp:e1}.

The proof of \eqref{app:l2:upp:e2} is similar to the proof of \eqref{app:l2:upp:e1}, but uses 
$ 1\leq 2^7\{|\widehat{w}_k/w_k-1|^8+ |\widehat{w}_k/w_k|^4|\widehat{w}_k/w_k-1|^4+|\widehat{w}_k/w_k|^4\}$ rather than \eqref{app:l2:upp:e1:1} and we omit the details.

Proof of \eqref{app:l2:upp:e3}. Due to the elementary inequality $1\leq 2|\widehat{c}_k/c_k-1|^2 + 2|\widehat{c}_k/c_k|^2$ we have
\begin{multline*}
\Ex\Bigl[\frac{|T_{n,k}|^2}{\widehat{w}_k\cdot \widehat{\lambda}_k} \1\{\widehat{w}_k\geq\alpha\}\1\{\widehat{\lambda}_k\geq \alpha\, \gamma_k^\nu\}\Bigr]\leq 
2 \Ex\Bigl[\frac{|T_{n,k}|^{2}|\widehat{c}_k/c_k-1|^2\1\{\widehat{w}_k\geq \alpha\}}{\widehat{w}_k\,\alpha\, \gamma_k^\nu} \Bigr]  + 2 \frac{\Ex|T_{n,k}|^{2}}{c_k^2}.
\end{multline*}
Consequently, \eqref{app:l2:upp:e1} and \eqref{app:l1:upp:e1} in Lemma \ref{app:l1:upp} lead to  \eqref{app:l2:upp:e3}.

Proof of \eqref{app:l2:upp:e4}. By using the Cauchy-Schwarz inequality we obtain the decomposition
\begin{multline*}
\Ex \Bigl[\Bigl|\frac{T_{n,k}\widehat{c}_k}{\widehat{w}_k}\Bigr|^2\1\{\widehat{w}_k\geq \alpha\}\1\{\widehat{\lambda}_k\geq \alpha\,\gamma_k^\nu\}\Bigr] \\
\leq 2\lambda_k 
 \Bigl(\Ex \frac{|T_{n,k}|^{4} w_k^2\1\{\widehat{w}_k\geq \alpha\}}{\widehat{w}_k^4}\Bigr)^{1/2} \cdot \Bigl\{
(\Ex|\widehat{c}_k/c_k-1|^4)^{1/2} + |P(\widehat{\lambda}_k\geq \alpha\,\gamma_k^\nu)|^{1/2}\Bigr\}.
\end{multline*}
Now \eqref{app:l1:upp:e3} in Lemma \ref{app:l1:upp},  \eqref{app:l3:upp:e3} in Lemma \ref{app:l3:upp} and \eqref{app:l2:upp:e2}  imply  \eqref{app:l2:upp:e4}, which completes the proof.\hfill$\square$\\ 
%%%%%%%%%%%%%%%%%%%%%%%%%%%%%%%%%%%%%%%%%
%%%%%%%%%%%%%%%%%%%%%%%%%%%%%%%%%%%%%%%%%%%
%%%%%%%%%%%%%%%%%%%%%%%%%%%%%%%%%%%%%%%%%%%
\begin{lem}\label{app:lem1:lower}Let $(\lambda_k)_{k\in \Z}$ be an element of $\cS_{\kappa,d}$ defined in \eqref{res:link:gen} with $\lambda_+:=1\vee\max_{k\in\Z}\lambda_k$. Consider $k^*\in\N$ and $\delta^*\in(0,1]$ given in \eqref{res:lower:def:md} for some $\triangle\geq 1$.  If we define \begin{equation}\label{app:lem1:lower:b}b_k^2:= \frac{\zeta}{n\cdot \lambda_k},\quad k\in\Z, \quad \text{ with }\quad \zeta:=\min \left\{ 2\sigma^2,  \rho/(d\, \triangle)\right\},\end{equation}
then we have
\begin{gather}\label{app:lem1:lower:e1}
\frac{n  }{2\sigma^2}b^2_k  \lambda_{k}\leq 1,\quad k\in\Z,\\\label{app:lem1:lower:e2}
\sum_{|k|\leq k^*}b^2_k   \gamma_{k}^{p}\leq \rho,\\\label{app:lem1:lower:e3}
 \sum_{|k|\leq k^*}b^2_k  \gamma_{k}^{\nu}\geq \min \left\{ \frac{2\sigma^2}{d\, \triangle},  \frac{\rho}{(d\, \triangle)^2}\right\} \cdot\frac{\kappa(\delta^*) }{\lambda_+}.
\end{gather}
\end{lem}
\noindent\textcolor{darkred}{\sc Proof.} The inequality \eqref{app:lem1:lower:e1} follows trivially by using the definition of $\zeta$.

Proof of \eqref{app:lem1:lower:e2}. If $\varphi$ denotes the inverse function of $\kappa$, then the link condition  $(\lambda_k)\in\cS_{\kappa,d}$, can  be rewritten as 
\begin{equation}\label{app:lem1:lower:e4}d^{-1}\leq |\varphi(\gamma_k^{\nu-p})|^{-1} \frac{\lambda_k}{\gamma_k^\nu\,\lambda_+}\leq d . \end{equation}
Thereby,  the monotonicity of $(\gamma_k^\nu)$  together with $\lambda_+\geq 1$ implies
\begin{equation*}\sum_{|k|\leq k^*}b^2_k \gamma_k^p \leq \frac{\zeta}{n}\cdot \sum_{|k|\leq k^*}\frac{\gamma_k^\nu\,\lambda_+}{\lambda_k}  \gamma_k^{p-v}
\leq \zeta\,d\cdot \sum_{|k|\leq k^*}\frac{\gamma_{k^*}^{p-v}}{n\varphi(\gamma_k^{\nu-p})} .
\end{equation*}
Thus \eqref{app:lem1:lower:e2} follows from the definition of $k^*$ given in \eqref{res:lower:def:md}, i.e.,  $\sum_{|k|\leq k^*}b^2_k \gamma_k^p \leq \zeta\,d\,\triangle\leq \rho$.

Proof of \eqref{app:lem1:lower:e3}. By using the condition \eqref{app:lem1:lower:e4} together with the definition of $\delta^*$ we have%
\begin{equation*} \sum_{|k|\leq k^*}b^2_k \gamma_k^\nu \geq \frac{\zeta}{d\,\lambda_+}\cdot\kappa(\delta^*)\cdot \sum_{|k|\leq k^*}\frac{\gamma_{k^*}^{p-v}}{n\varphi(\gamma_k^{\nu-p})}.
\end{equation*}
Consequently, the definition of $k^*$ given in \eqref{res:lower:def:md} implies \eqref{app:lem1:lower:e3}, which proves the lemma.\hfill$\square$%\\

%% file: _a5proofs-infin.tex
\paragraph{The finitely smoothing case.}\hfill\\[1ex]
\noindent\textcolor{darkred}{\sc Proof of Proposition \ref{res:lower:pol}.} Consider the inverse function $\varphi(t)=t^{(a+\nu)/(p-\nu)}$ of $\kappa$. Then  the well-known  approximation $\sum_{k=1}^{m} k^{r}\asymp m^{r+1}$ for $r>0$ together with the definition of $\gamma_k$ given in \eqref{bm:def:gamma} implies 
$\sum_{|k|\leq k^*}1/\varphi(\gamma_k^{\nu-p})\asymp \gamma_{k^*}^{(a+\nu)+1/2}$. It follows that the condition on $k^*$  given in \eqref{res:lower:def:md} of Theorem \ref{res:lower} can be rewritten as 
\begin{equation*}
1/n\asymp \gamma_{k^*}^{\nu-p}\Bigl|\sum_{|k|\leq k^*}1/\varphi(\gamma_k^{\nu-p})\Bigr|^{-1}\asymp \gamma_{k^*}^{p+a+1/2}=|\varphi 
(\gamma_{k^*}^{\nu-p})|^{[2(p+a)+1]/[2(a+\nu)]} .
\end{equation*} 
From this   $\delta^*:=\varphi 
(\gamma_{k^*}^{\nu-p})$ implies $\delta^*\asymp n^{-2(a+\nu)/ [2(p+a)+1]}$ and  $\kappa(\delta^*)\asymp n^{-2(p-\nu)/ [2(p+a)+1]}$.  Consequently, the lower bound in Proposition \ref{res:lower:pol} follows by applying Theorem \ref{res:lower}. \hfill$\square$\\[1ex]
% %%%%%%%%%%%%%%%%%%%%%%%%%%%%%%%%%%%%%%%%%%
\noindent\textcolor{darkred}{\sc Proof of Proposition \ref{res:upper:pol}.} Since the condition on $\alpha$ ensures $\alpha \asymp n^{-2(a+\nu)/ [2(p+a)+1]} \asymp \delta^*$ (see the proof of Proposition \ref{res:lower:pol}) the result follows from Theorem \ref{res:upper}.
 \hfill$\square$%\\
% %%%%%%%%%%%%%%%%%%%%%%%%%%%%%%%%%%%%%%%%%%
\paragraph{The infinitely smoothing case.}\hfill\\[1ex]
\noindent\textcolor{darkred}{\sc Proof of Proposition \ref{res:lower:log}.} Consider the inverse function $\varphi(t)=\exp\{-t^{a/(\nu-p)}\}$ of $\kappa$. By applying Laplace's Method (c.f. chapter 3.7 in \cite{Olver1974}) the definition of $\gamma_k$ given in \eqref{bm:def:gamma} implies
$\sum_{|k|\leq k^*}1/\varphi(\gamma_k^{\nu-p})\asymp 1/\varphi(\gamma_{k^*}^{\nu-p})$. It follows that by using $\delta^*:=\varphi 
(\gamma_{k^*}^{\nu-p})$ the condition on $k^*$  given in \eqref{res:lower:def:md} of Theorem \ref{res:lower} can be rewritten as 
\begin{equation*}
1/n\asymp \gamma_{k^*}^{\nu-p}\Bigl|\sum_{|k|\leq k^*}1/\varphi(\gamma_k^{\nu-p})\Bigr|^{-1}\asymp \gamma_{k^*}^{\nu-p} \varphi(\gamma_{k^*}^{\nu -p})=\delta^*\kappa(\delta^*),\end{equation*} 
which implies $\kappa(\delta^*)\asymp \omega(1/n)$, where $\omega$ denotes the inverse function of $\omega^{-1}(t)=t\cdot\varphi(t)$. 
Therefore, the lower bound given in Proposition \ref{res:lower:log} follows  from Theorem \ref{res:lower} together with 
$\omega(t)=|\log t|^{-(p-\nu)/a}(1+o(1))$ as $t \to0$ (c.f. \cite{Mair94}), which  proofs the result \hfill$\square$\\[1ex]
% %%%%%%%%%%%%%%%%%%%%%%%%%%%%%%%%%%%%%%%%%%
\noindent\textcolor{darkred}{\sc Proof of Proposition \ref{res:upper:log}.} The proof is based on  the  decomposition \eqref{res:dec}, where we bound the two right hand side terms  by \eqref{pr:res:prop:cons:e1} derived in the proof of Proposition \ref{res:prop:cons} and \eqref{pr:res:upper:gen:e2} shown  in  the proof of Theorem \ref{res:upper} respectively. It follows that,
\begin{equation*}\Ex\|\widehat{\beta}_{\nu}-\beta\|_{\nu}^2\leq \frac{C\,\eta}{\alpha^2 \, n}\cdot \Lambda \cdot\{\sigma^2+\rho\Lambda\}  + \frac{C\,\eta }{\alpha\, n}\cdot\Lambda\cdot \rho \cdot\Bigl(1+\frac{\Lambda}{\alpha\, n}\Bigr) + C\cdot \kappa(d\,4\,\tau\, \alpha)\cdot \rho
\end{equation*}
 for some positive constant $C$. Consequently,  the condition  $\alpha=c\cdot n^{-1/4},$ $c>0$ implies 
$ \Ex\|\widehat{\beta}_{\nu}-\beta\|_{\nu}^2 = O( n^{-1/2})  + O(n^{-3/4}) + O( |\log n^{-1/4}|^{-(p-\nu)/a})=O((\log n)^{-(p-\nu)/a}))$, which completes the proof. \hfill$\square$%\\
% %%%%%%%%%%%%%%%%%%%%%%%%%%%%%%%%%%%%%%%%%%